\documentclass[a4paper,11pt]{article}
\usepackage[latin1]{inputenc}
\usepackage[T1]{fontenc}
\usepackage{a4wide}
\usepackage{fontenc}
\usepackage{color}
\usepackage{amsmath}
\usepackage{amsthm} 
\usepackage{graphicx}
\usepackage{array}
\usepackage{amsfonts}
\usepackage{amssymb}
\usepackage{euscript}
\usepackage{delarray}
\usepackage{latexsym}
\usepackage{enumerate}
\usepackage{float}
\usepackage{algorithmicx}
\usepackage{algorithm}
\usepackage{algpseudocode}
\usepackage{url}
\usepackage{caption}
\usepackage{subcaption}
\usepackage{textcomp}

\title{A Bayesian approach for global sensitivity analysis  of (multi-fidelity) computer codes}
\author{Loic Le Gratiet $^\dag$  $^\ddag$, Claire Cannamela $^\ddag$, Bertrand Iooss $^\S$ $^\star$ \\ \\ $^\dag$ Universit\'e Paris Diderot 75205 Paris Cedex 13, France \\ \\ $^\ddag$ CEA, DAM, DIF, F-91297 Arpajon, France \\ \\ $^\S$ EDF R\&D, 6 quai Watier 78401 Chatou, France \\ \\ $^\star$ Institut de Math\'ematiques de Toulouse, 31062 Toulouse, France  }

\newtheorem{prop}{Proposition}

\newcommand{\EX}[1]{\mathbb{E}_X\left[#1\right]}
\newcommand{\VX}[1]{\mathrm{var}_X\left(#1\right)}
\newcommand{\covX}[2]{\mathrm{cov}_X\left(#1,#2\right)}
\newcommand{\EZ}[1]{\mathbb{E}_Z\left[#1\right]}
\newcommand{\VZ}[1]{\mathrm{var}_Z\left(#1\right)}

\newcommand{\N}[2]{\mathcal{N}\left(#1,#2\right)}

\newcommand{\eqL}[0]{\stackrel{\mathcal{L}}{=}}
\newcommand{\tx}[0]{{\tilde{x}}}
\newcommand{\tX}[0]{{\tilde{X}}}
\newcommand{\tildt}[0]{{\tilde{t}}}

\newcommand{\zn}[0]{\mathbf{z}^n}

\newcommand{\F}[0]{\mathbf{F}}

\newcommand{\D}[0]{\mathbf{D}}

\newcommand{\xx}[0]{\mathbf{x}}

\newcommand{\rr}[0]{\mathbf{r}}

\newcommand{\hh}[0]{\mathbf{h}}
\newcommand{\ff}[0]{\mathbf{f}}

\newcommand{\zz}[0]{\mathbf{z}}

\newcommand{\RR}[0]{\mathbf{R}}
\newcommand{\ZZ}[0]{\mathbf{Z}}

\newcommand{\HHH}[0]{\mathbf{H}}

\newcommand{\R}[0]{\mathbb{R}}

\newcommand{\rrho}[0]{\boldsymbol{\rho}}

\newcommand{\bbeta}[0]{\boldsymbol{\beta}}

\newcommand{\ssigma}[0]{\boldsymbol{\sigma}}

\newcommand{\ttheta}[0]{\boldsymbol{\theta}}

\begin{document}
\maketitle

		\section*{Abstract}

		Complex computer codes  are widely used in science and engineering to model physical phenomena. Furthermore,  it is common that they have a large number of input parameters. Global sensitivity analysis aims to identify those which have the most important impact on the output.
		Sobol indices are  a popular    tool   to perform such analysis. However, their  estimations require  an important number of simulations and often cannot be processed under reasonable time constraint.
		To handle this problem, a Gaussian process regression model is built  to surrogate the computer code and the Sobol indices are estimated through it. The  aim of this paper is to provide a methodology to estimate the Sobol indices through a surrogate model  taking into account both the estimation errors and the surrogate model errors. 
		In particular, it allows us to derive non-asymptotic confidence intervals for the Sobol index estimations.
		Furthermore,  we extend the suggested strategy   to the case of  multi-fidelity computer codes which can be run at different levels of accuracy.
		For such simulators, we use an  extension of Gaussian process regression models  for multivariate outputs.

		\paragraph{Keywords:} Sensitivity analysis, Gaussian process regression, multi-fidelity model, complex computer codes, Sobol index, Bayesian analysis.

		\section{Introduction}

		Complex computer codes commonly have a large number of input parameters for which we want to measure their importance on the model response. We focus on the  Sobol indices   \cite{sobol1993}, \cite{Sal00} and  \cite{sobol2}   which are  variance-based importance measures  coming from the Hoeffding-Sobol decomposition \cite{hoeffding1948}. We note that  the presented sensitivity analysis holds when the input parameters are independent. For an analysis with dependent inputs, the reader is referred to the articles \cite{kucherenko}  and  \cite{chastaing}.

		A widely used method to estimate the Sobol indices are the Monte-Carlo based methods.   They allow  for quantifying the errors due to numerical integrations (with a Bootstrap procedure in a non-asymptotic case \cite{archer1997sensitivity} and  \cite{janon2011uncertainties} or thanks to asymptotic normality properties in an asymptotic case \cite{Jan12}).
		However, the estimation of the Sobol indices by sampling methods requires a large number of simulations that are sometimes too costly and time-consuming. A popular method to overcome this difficulty is to build a mathematical approximation of the code input/output relation \cite{marseguerra2003variance} and \cite{iooss2006response}. 

		We deal  in this paper with the  use of kriging and multi-fidelity co-kriging models  to  surrogate the computer code. The reader is referred to the books \cite{S99}, \cite{S03} and \cite{R06} for an overview of kriging methods for computer experiments.  A pioneering article dealing with the kriging approach to perform global sensitivity analysis is the one of Oakley and O'Hagan \cite{Oak04}. Their  method is also investigated in \cite{Mar09}.    The strength of the suggested approach is that it allows for inferring from the surrogate model uncertainty about the Sobol index estimations. 
		However, it does not use Monte-Carlo integrations and   it does not take into account the numerical errors due to the   numerical integrations. 
		Furthermore, the    implementation of the method is complex  for general covariance kernels. 
		Another flaw of the method presented in \cite{Oak04} and \cite{Mar09} is that it is not based on the exact definition of Sobol indices (it uses the ratio of two expectations instead of the expectation of a ratio).%able to handle the real Sobol  indices.
		We note that a bootstrap procedure can also be used to evaluate the impact of the  surrogate model uncertainty on the Sobol index estimates as presented in \cite{storlie2009implementation}. However, this approach only  focuses on the parameter estimation errors.

		On the other hand, a method   giving confidence intervals for the Sobol  index estimations and taking into account  both   the meta-model uncertainty and the numerical integration errors   is suggested in  \cite{janon2011uncertainties}.  
They consider   Monte-Carlo integrations   to estimate the Sobol  indices (see  \cite{sobol1993},  \cite{sobol2007estimating} and \cite{Jan12})  instead of numerical integrations   and they infer from the sampling errors thanks to  a bootstrap procedure. Furthermore, to deal with the meta-model error, they consider an upper bound on it.  In the kriging case  they use the kriging variance up to a  multiplicative constant as upper bound.  Nevertheless, this is a rough upper bound which considers the worst error on a test sample.
		Furthermore,  this method does not  allow for inferring from the meta-model uncertainty about the Sobol   index estimations.

		We propose in this paper a method combining the approaches  presented in  \cite{Oak04} and  \cite{janon2011uncertainties}. As in \cite{Oak04} we consider  the code as  a realization of a Gaussian process. Furthermore,  we use the method  suggested in \cite{janon2011uncertainties} to estimate the Sobol  indices with Monte-Carlo integrations.  Therefore, we can use the bootstrap  method presented in  \cite{archer1997sensitivity}  to infer from the sampling error  on    the Sobol   index estimations. Furthermore, contrary to \cite{Oak04} and \cite{Mar09} we deal with the exact definition of Sobol   indices. Consequently, we introduce   non-asymptotic  certified Sobol  index  estimations, i.e. with  confidence intervals which take into account both the surrogate model error and the numerical integration errors.

		Finally, the suggested approach is extended  to a multi-fidelity   co-kriging model. 
		It allows for approximating a computer code using fast and coarse versions of it. The suggested multi-fidelity models is derived from the original one  proposed in \cite{KO00}. We note that the use of co-kriging model to deal with multi-fidelity   codes have been largely investigated during this  last decade  (see \cite{KO01}, \cite{Hig04}, \cite{Ree04} and  \cite{QW07}).
		 A  definition of Sobol   indices  for multi-fidelity computer codes is  presented in \cite{Jac06}. However, their approach  is based on  tabulated  biases between fine and coarse codes and does not allow for inferring from the meta-model uncertainty. The co-kriging model fixes these  weakness since it allows for considering general forms for the biases and for inferring from the surrogate model error.

		This  paper is organized as follows. First we introduce in Section \ref{chap6sec1} the so-called Sobol indices. Then, we present in Section \ref{chap6sec2}   the   kriging-based sensitivity analysis suggested  by  \cite{Oak04}. Our approach is developed in Section \ref{chap6sec3}. In particular, we give an important result allowing for effectively sampling with respect to the kriging predictive distribution in Subsection \ref{chap6sec30}. Finally, we extend in Section \ref{chap6sec4} the presented approaches to multi-fidelity co-kriging models. We highlight that we present in  Subsection \ref{chap6sec40} a method to sampling with respect to the multi-fidelity predictive distribution. In this case the predictive distribution is not anymore Gaussian.  Numerical tests are performed in Section \ref{sec:tests} and an industrial example is considered in Section \ref{sec:application}. A conclusion synthesizes this work in the last section.
		
		%%%%%%%%%%%%%%%%%%%%%%%%
		\section{Global sensitivity analysis: the method of Sobol}\label{chap6sec1}

		We present in this section the method of Sobol for global sensitivity analysis \cite{sobol1993}. It is inspired by the book of \cite{Sal00} giving an   overview of classical sensitivity analysis methods.

			\subsection{Sobol variance-based sensitivity analysis}

			Let us  consider the input parameter space  $Q \subset \R^d$ such  that $(Q,\mathcal{B}(Q))$ is a  measurable product space of the form:
			\begin{displaymath}
			(Q,\mathcal{B}(Q)) = (Q_1 \times  \dots \times Q_d, \mathcal{B}(Q_1 \times  \dots \times Q_d))
			\end{displaymath}
			where $\mathcal{B}$ is the Borelian $\sigma$-algebra and $Q_i \subset \R$ is a nonempty open set, for $i=1,\dots,d$. Furthermore,  we consider a probability measure $ \mu  $ on  $(Q,\mathcal{B}(Q))$, with values in $\R$ and of the form
			\begin{displaymath}
			\mu(x) = \mu_1(x^1) \otimes \dots \otimes \mu_d(x^d)
			\end{displaymath}
			The Hoeffding-Sobol  decomposition (see \cite{hoeffding1948})  states that any function $z(x) \in L^2_\mu(\R^d)$ can   be decomposed into summands of increasing dimensionality in  such way:
			\begin{equation}\label{sobolDec}
			z(x) = z_0 + \sum_{i=1}^d z_i(x^i) + \sum_{1\leq i < j \leq k}z_{ij}(x^i,x^j)+\dots+z_{1,2,\dots,d}(x^1,\dots,x^d) = \sum_{u \in  \mathcal{P}} z_u(x^u)
			\end{equation}
			where $\mathcal{P}$ is  the collection of all subsets of $ \{1,\dots,d\}$ and  $x^u$ is a group of variables such that  $x^u = (x^i)_{i \in u}$. Furthermore, the decomposition is unique if we consider  the following  property for every summand $u = (u_1,\dots,u_k)_{\substack{1\leq k \leq d }}$, $1 \leq u_i \leq d$:
			\begin{equation}
			\int z_u(x^u)  \,d\mu_{u_i}(x^{u_i}) = 0, \quad \forall i=1,\dots,k.
			\end{equation}
%			A consequence of this property is that all the summands are orthogonal, i.e.  for every   $z_u(x^u)$ and $z_v(x^v)$ such that  $u, v \in \mathcal{P}$ and $u \neq v$, we have:
%			\begin{equation}\label{orthoSobol}
%			\int z_u(x^u) z_v(x^v) \,d\mu(x) = 0
%			\end{equation}
%			Another consequence is that $z_0$ represents the mean of $z(x)$ with respect to the measure $\mu(x)$
%			\begin{equation}
%			z_0 = \int z(x) \, d\mu(x)
%			\end{equation}	
			Now, let us suppose that  the inputs  are  a random vector $X = (X^1,\dots,X^d)$  defined on the probability space $(\Omega_X, \mathcal{F}_X, \mathbb{P}_X)$  and with  measure $\mu$. Sobol \cite{sobol1993} showed that the decomposition (\ref{sobolDec})  can be interpreted as conditional expectations as follows:
			\begin{eqnarray*}
			z_0 & = & \EX{z(X)}\\
			z_i(X^i) & = &  \EX{z(X) |X^i} - z_0 \\
			z_{ij}(X^i,X^j) & = & \EX{z(X) |X^i,X^j} - z_i(X^i) - z_j(X^j) - z_0 \\
			 & \vdots & \\
			z_u(X^u) & = & \EX{z(X) |X^u} - \sum_{v \subset u} z_v(X^v)
			\end{eqnarray*}		
			with $u  \in  \mathcal{P}$. From this scheme, we can naturally develop the variance-based sensitivity indices of Sobol. First, let us consider the total variance $D$ of $z(x)$:
			\begin{equation}\label{totalvariance}
			D = \VX{z(X)}
			\end{equation} 
			By squaring and integrating the decomposition (\ref{sobolDec}), we obtain 
			\begin{equation}
			D = \sum_{i=1}^d D_i + \sum_{1 \leq i< j \leq d } D_{ij} + \dots + D_{1,2,\dots,d} = \sum_{u \in \mathcal{P}} D_u.
			\end{equation}
			with 
			$
			D_u = \VX{\EX{z(X) |X^u}} - \sum_{v \subset u}\VX{\EX{z(X) |X^v}}  .
			$
			Finally, the Sobol sensitivity indices are given by
			\begin{equation}
			S_u = \frac{D_u}{D}
			\end{equation}
			where $u \in \mathcal{P}$. 
			We note that we have the following useful equality which allows for easily interpreting $S_u$ as the part of variance  of $z(x)$  due to  $x^u$ and not explained by $x^v$ with $v \subset u$.
			\begin{equation}
			1 = \sum_{i=1}^d S_i + \sum_{1 \leq i< j \leq d } S_{ij} + \dots + S_{1,2,\dots,d} = \sum_{u \in \mathcal{P}} S_u.
			\end{equation}
			In particular, $S_i$ is called the first-order sensitivity index for variable $x^i$. It measures the main effect of $x^i$ on the output, i.e. the part of variance  of $z(x)$ explained by the factor $x^i$. Furthermore, $S_{ij}$ for $i \neq j$ is the second-order sensitivity index. It measures the part of variance  of  $z(x)$ due to $x^i$ and $x^j$  and not explained by the individual effects of $x^i$ and $x^j$.  

			\subsection{ Monte-Carlo Based estimations of Sobol indices}\label{MCsobolestim}

			Now, let us denote by $Q^{d_1} = Q_{i_1} \times  \dots \times Q_{i_{d_1}}$, $d_1 \leq d$, $\{i_1,\dots,i_{d_1}\} \in \mathcal{P}$ and $Q^{d_2} = Q_{j_1} \times  \dots \times Q_{j_{d_2}}$ such that $\{j_1,\dots,j_{d_2}\} = \{1,\dots,d\} \setminus \{i_1,\dots,i_{d_1}\}$. Analogously, we  use the notation $X^{d_1} = (X^i)_{i \in \{i_1,\dots,i_{d_1}\}}$, $X^{d_2} = (X^j)_{j \in \{j_1,\dots,j_{d_2}\}}$,  $\mu^{d_1} = \left( \bigotimes_{i \in \{i_1,\dots,i_{d_1}\}} \mu_i \right) $ and $\mu^{d_2} =  \left(\bigotimes_{j \in \{ j_1,\dots,j_{d_2}\}} \mu_j \right)$ where  $\mu^{d_1}$ and  $\mu^{d_2}$ are   probability measures on $(Q^{d_1}, \mathcal{B}(Q^{d_1}))$ and $(Q^{d_2}, \mathcal{B}(Q^{d_2}))$. Consequently, we have the equalities   $\mu = \mu^{d_1}
 \otimes \mu^{d_2}$,  $Q = Q^{d_1} \times Q^{d_2}$ and $X = (X^{d_1}, X^{d_2})$ with $d = d_1 + d_2$. 

			We are interested in evaluating the closed sensitivity index:
			\begin{equation}\label{VX1V}
			\mathcal{S}^{X^{d_1}} = \frac{V^{X^{d_1}}}{V} = \frac{\VX{\EX{z(X)|X^{d_1}}}}{\VX{z(X)}}
			\end{equation}
			A first method would be to use $d$-dimensional numerical integrations to approximate the numerator and denominator of (\ref{VX1V}) as presented in \cite{Oak04} and \cite{Mar09}. Nonetheless, since $d$ is large in general, this method leads to  numerical issues and is computationally expensive. A second approach is to take advantage of the probabilistic interpretation of the Sobol indices and to use a Monte-Carlo procedure to evaluate the different integrals as presented in the forthcoming developments.  

			\begin{prop}\label{varcovsob}
			Let us consider  the random vectors $(X, \tX)$ with $X = (X^{d_1}, X^{d_2})$  and $\tX = (X^{d_1}, \tX^{d_2})$ where $X^{d_1}$ is a random vector  with   measure $\mu^{d_1}$ on $Q^{d_1}$, $X^{d_2}$ and $\tX^{d_2}$ are  random vectors   with   measure $\mu^{d_2}$  on $Q^{d_2}$  and $X^{d_2}  \perp  \tX^{d_2}$. We have the following equality:
			\begin{equation}\label{covX1}
			\VX{\EX{z(X)|X^{d_1}}} = \covX{z(X)}{z(\tX)}
			\end{equation}
			\end{prop}

			$\mathcal{S}^{X^{d_1}} $ in equation (\ref{VX1V}) can thus be estimated by considering two  random vectors  $(X_i)_{i=1,\dots,m}$ and $(\tX_i)_{i=1,\dots,m}$, $m \in \mathbb{N}^*$ lying  in $(\Omega_X, \mathcal{F}_X, \mathbb{P}_X)$ such that  $X_i \eqL X$ and $\tX_i \eqL \tX$ ($\eqL$ stands for an equality in distribution) and by  using an estimator for the covariance $\covX{z(X)}{z(\tX)}$.

			Following this principle, Sobol \cite{sobol1993} suggests the following estimator for the ratio  in equation (\ref{VX1V}):
			\begin{equation}\label{Sobolestim}
			\frac{V_m^{X^{d_1}}}{V_m} = \frac{\frac{1}{m} \sum_{i=1}^{m} z(X_i)z(\tX_i)
			-
			 \frac{1}{m}\sum_{i=1}^{m} z(X_i)  \frac{1}{m}\sum_{i=1}^{m}z(\tX_i) 
			}{
			\frac{1}{m} \sum_{i=1}^{m} z(X_i)^2
			-
			\left( \frac{1}{m}\sum_{i=1}^{m} z(X_i)\right)^2
			} 
			\end{equation}
			This estimation is improved by  \cite{Jan12} who  propose the following estimator:
			\begin{equation}\label{Janonestim}
			\frac{V_m^{X^{d_1}}}{V_m} = \frac{\frac{1}{m} \sum_{i=1}^{m} z(X_i)z(\tX_i)
			-
			\left(  \frac{1}{2m}\sum_{i=1}^{m} z(X_i) + z(\tX_i)   \right)^2  
			}{
			\frac{1}{m} \sum_{i=1}^{m} z(X_i)^2
			-
			\left(  \frac{1}{2m}\sum_{i=1}^{m} z(X_i) + z(\tX_i)   \right)^2  
			} 
			\end{equation}
			In particular they demonstrate that  the asymptotic variance in (\ref{Janonestim})  is better than the one  in (\ref{Sobolestim}) and they show that the estimator (\ref{Janonestim})  is asymptotically efficient for the first order indices. The main weakness of the estimators   (\ref{Sobolestim}) and (\ref{Janonestim})  is that they are sometimes not accurate for small values of ${V^{X^{d_1}}}/{V}$ in (\ref{VX1V}). To tackle this issue, \cite{sobol2007estimating} propose the following estimator
			\begin{equation}\label{Mauntzestim}
			\frac{V_m^{X^{d_1}}}{V_m} = \frac{\frac{1}{m} \sum_{i=1}^{m} z(X_i)z(\tX_i)
			-
			 \frac{1}{m}\sum_{i=1}^{m} z(X_i) z(\tilde{\tX}_i) 
			}{
			\frac{1}{m} \sum_{i=1}^{m} z(X_i)^2
			-
			\left( \frac{1}{m}\sum_{i=1}^{m} z(X_i)\right)^2
			} 
			\end{equation}
			where $\tilde{\tX} = (\tX^{d_1}, \tX^{d_2})$, $\tX^{d_1} \eqL X^{d_1}$, $\tX^{d_1} \perp  X^{d_1}$  and $(\tilde{\tX}_i)_{i=1,\dots,m}$   is such that  $\tilde{\tX}_i  \eqL \tilde{\tX}$ for all $i=1,\dots,m$.

		\section{Kriging-based sensitivity analysis: a first approach}\label{chap6sec2}

		We present in this Section the approach  suggested in \cite{Oak04}  and \cite{Mar09} to perform global sensitivity analysis using kriging surrogate models. Then, we present an alternative  method   that  allows for avoiding  complex numerical integrations.  Nevertheless, we will see that this approach does not provide a correct representations  of the Sobol indices. We handle this problem in the next section.

			\subsection{A short introduction to kriging model}

			The principle of the kriging model   is to consider that our prior knowledge about  the code $z(x)$  can be modelled by a Gaussian process $Z(x)$ with mean $\ff'(x) \bbeta$ and covariance kernel $\sigma^2 r(x,\tx)$ (see for example \cite{S03}). Then,   the code $z(x)$ is approximated by a Gaussian process $Z_n(x)$ having the predictive distribution $[ Z(x)|Z(\D) =  \zn , \sigma^2]$   where $\zn$ are the known values of  $z(x)$ at points in  the experimental design set $\D = \{ x^1,\dots,x^n\}$, $x^i \in Q$, and   $\sigma^2$ is the variance parameter:
			\begin{equation}\label{predictiveZn}
			Z_n(x) \sim \mathrm{GP}\left(m_n(x), s_n^2(x,\tx) \right)
			\end{equation}
			where the mean $m_n(x)$ and the variance $s_n^2(x,\tx)$  are given by:
			\begin{displaymath}
			m_n(x) = \ff'(x) \hat{\bbeta} + \rr'(x) \RR^{-1}\left( \zn - \F \hat{\bbeta} \right)
			\end{displaymath}
			where $\RR=[r(x_i,x_j)]_{i,j=1,\dots,n}$, $\rr'(x)=[r(x,x_i)]_{i=1,\dots,n}$, $\F =[\ff'(x_i)]_{i=1,\dots,n}$ and 
			\begin{displaymath}
			s^2_n(x,\tx) =\sigma^2\left( 1 - \begin{pmatrix} \ff'(x) & \rr'(x)  \end{pmatrix} \begin{pmatrix} 0 & \F' \\ \F & \RR \end{pmatrix}^{-1} \begin{pmatrix} \ff(\tx) \\ \rr(\tx) \end{pmatrix}\right)
			\end{displaymath}
			where $\hat{\bbeta} = \left( \F' \RR^{-1} \F \right)^{-1} \F' \RR^{-1} \zn$. The variance parameter $\sigma^2$ can be estimated with a restricted maximum likelihood method, i.e. $\hat{\sigma}^2 = (\zn - \hat{\bbeta}\F)' \RR^{-1}(\zn - \hat{\bbeta}\F)/(n-p)$ where $p$ is the size of $\bbeta$.

			\subsection{Kriging-based Sobol index}\label{oakleymethod}

			The idea suggested in \cite{Oak04} and \cite{Mar09} is to substitute $z(x)$ with $Z_n(x)$ in equation (\ref{VX1V}): 
			\begin{equation}\label{VX1VZ}
				\mathcal{S}_n^{X^{d_1}} = \frac{V_n^{X^{d_1}}}{V_n} = \frac{\VX{\EX{Z_n(X)|X^{d_1}}}}{\VX{Z_n(X)}}
			\end{equation}
			Therefore, if we denote by $(\Omega_Z, \mathcal{F}_Z, \mathbb{P}_Z)$ the probability space where the Gaussian process $Z(x)$ lies, then the estimator $\mathcal{S}_n^{X^{d_1}}$ lies in $(\Omega_Z, \mathcal{F}_Z, \mathbb{P}_Z)$ (it is hence random). We note that  $Z_n(X)$  is defined on the product probability space $(\Omega_X \times \Omega_Z, \sigma(\mathcal{F}_X \times \mathcal{F}_Z), \mathbb{P}_X \otimes \mathbb{P}_Z)$. 

			Nevertheless,  the distribution of $\mathcal{S}_n^{X^{d_1}}$ is intractable and \cite{Oak04} and \cite{Mar09}   focus on  its mean and variance. More precisely, in order to derive analytically the Sobol index  estimations they consider the following quantity:
			\begin{equation}\label{meanoak}
			\tilde{\mathcal{S}}_n^{X^{d_1}} = \frac{\EZ{\VX{\EX{Z_n(X)|X^{d_1}}}}}{\EZ{\VX{Z_n(X)}}}
			\end{equation}
			where $\EZ{.}$ stands for the expectation in the probability space  $(\Omega_Z, \mathcal{F}_Z, \mathbb{P}_Z)$. Furthermore, the uncertainty on $\tilde{\mathcal{S}}_n^{X^{d_1}}$ is  evaluated  with the following quantity:
			\begin{equation}\label{varoak}
			\sigma^2( \tilde{\mathcal{S}}_n^{X^{d_1}}) = \frac{\VZ{\VX{\EX{Z_n(X)|X^{d_1}}}}}{\EZ{\VX{Z_n(X)}}^2}
			\end{equation}
			As shown  in \cite{Oak04} and \cite{Mar09}, the equations (\ref{meanoak}) and (\ref{varoak}) can be derived analytically through multi-dimensional integrals for the cases $d_1 = i$, $i=1,\dots,d$, i.e. for the first-order indices. Furthermore, with some    particular  formulations  of  $\ff(x)$, $\mu(x)$ and $r(x,\tx)$,  these multi-dimensional  integrals can be written as product of one-dimensional ones.  We note that a method is suggested in \cite{Mar09} to generate samples  of the numerator $\VX{\EX{Z_n(X)|X^{d_1}}}$ in (\ref{VX1VZ}). It allows for estimating    the uncertainty of   $\tilde{\mathcal{S}}_n^{X^{d_1}}$ in  (\ref{meanoak}) without processing the complex numerical integrations involved in (\ref{varoak}).

			\paragraph{Discussions:} The method suggested in \cite{Oak04} and \cite{Mar09} provides an interesting tool to perform sensitivity analysis of complex models. Nevertheless, in our opinion it suffers from the following flaws:
			\begin{enumerate}
			\item For general choice of $\ff(x)$, $\mu(x)$ and $r(x,\tx)$, the numerical evaluations of (\ref{meanoak}) and (\ref{varoak}) can be very complex   since it requires multi-dimensional integrals.
			\item The method  is derived for first-order sensitivity indices and cannot easily be  extended  to higher order indices.
			\item The method allows  for inferring from the surrogate model uncertainty  about the sensitivity indices   but does not allow for taking into account the numerical errors due to the multi-dimensional integral estimations.
			\item The considered index expectation and deviation  do  not correspond to the real Sobol index ones since we obviously have
			\begin{displaymath}
			  \frac{\EZ{\VX{\EX{Z_n(X)|X^{d_1}}}}}{\EZ{\VX{Z_n(X)}}} \neq   \EZ{ \frac{\VX{\EX{Z_n(X)|X^{d_1}} }}{ \VX{Z_n(X)}}}		
			\end{displaymath}
			and 
			\begin{displaymath}
			  \frac{\VZ{\VX{\EX{Z_n(X)|X^{d_1}}}}}{\EZ{\VX{Z_n(X)}}^2} \neq   \VZ{ \frac{\VX{\EX{Z_n(X)|X^{d_1}} }}{ \VX{Z_n(X)}}}		
			\end{displaymath}
			\end{enumerate}
		
			In the next subsection, we deal with the points 1, 2 and 3  by suggesting a Monte-Carlo sampling method to evaluate (\ref{meanoak}) and  (\ref{varoak}) instead of quadrature integrations. Nonetheless, we do not tackle the   issue of point 4. To handle it, we suggest another method in Section \ref{chap6sec3}.

			\subsection{Monte-Carlo estimations for the first approach}\label{MCoak}

			We present in this Subsection, another approach to   deal with the evaluation of $\tilde{\mathcal{S}}_n^{X^{d_1}} $ in (\ref{meanoak}). Its principle  simply consists in using the  estimation methods suggested in Subsection \ref{MCsobolestim} instead of quadrature integrations to compute ${\EZ{\VX{\EX{Z_n(X)|X^{d_1}}}}}$ and ${\EZ{\VX{Z_n(X)}}}$. We present the method with the estimator presented in \cite{sobol1993}. The extension to those presented in  \cite{janon2011uncertainties} and  \cite{sobol2007estimating} is straightforward. Let us substitute in  the estimator presented in equation  (\ref{Sobolestim}) the code $z(x)$ by  the Gaussian process $Z_n(x)$:
			\begin{equation}\label{SobolestimZ}
			\frac{V_{m,n}^{X^{d_1}}}{V_{m,n}} = \frac{\frac{1}{m} \sum_{i=1}^{m} Z_n(X_i)Z_n(\tX_i)
			-
			 \frac{1}{m}\sum_{i=1}^{m} Z_n(X_i)  \frac{1}{m}\sum_{i=1}^{m}Z_n(\tX_i) 
			}{
			\frac{1}{m} \sum_{i=1}^{m} Z_n(X_i)^2
			-
			\left( \frac{1}{m}\sum_{i=1}^{m} Z_n(X_i)\right)^2
			} 
			\end{equation}
			where the samples  $(X_i)_{i=1,\dots,m}$ and $(\tX_i)_{i=1,\dots,m}$  are those introduced in Subsection \ref{MCsobolestim}. Therefore, ${V_{m,n}^{X^{d_1}}}/{V_{m,n}}$ is an estimator of ${V^{X^{d_1}}}/{V}$ (\ref{VX1V})  when we replace the true function $z(x)$ by its approximation $Z_n(x)$ built from $n$ observations $\zn$ of $z(x)$ and when we estimate the variances and  the expectation  involved in  (\ref{VX1V}) by a Monte-Carlo method with $m$ particles.
To be clear in the remainder of this paper, we name  as Monte-Carlo error the one due to the Monte-Carlo estimation   and we name as  meta-model error the one due to the substitution of $z(x)$ by a surrogate model. Furthermore, $m$ will always denote the number of Monte-Carlo particles and $n$ the number of observations used to build the surrogate model.

			The strength of this formulation is that it gives closed form formulas for the evaluation of (\ref{meanoak})  for any choice of $\ff(x)$, $\mu(x)$ and $r(x,\tx)$ contrary to \cite{Oak04} and \cite{Mar09}. Furthermore, this method can directly be used for any order of Sobol indices which contrasts  with the one presented in Subsection (\ref{oakleymethod}). Finally, unlike quadrature integrations, Monte-Carlo integrations    allow  for taking into account the numerical errors due to the integral evaluations. In particular, as presented in \cite{archer1997sensitivity}, the bootstrap method can be directly used  to obtain confidence intervals on the Sobol indices.
	
			We give in the following equation the Monte-Carlo estimation  of $\tilde{\mathcal{S}}_n^{X^{d_1}} $  (\ref{meanoak}) corresponding to the  kriging-based sensitivity indices  presented in \cite{Oak04} and \cite{Mar09}.

			\begin{equation} \label{meanoakestim}
			\begin{array}{lll}
			\tilde{\mathcal{S}}_{m,n}^{X^{d_1}}  &  = & 
			\displaystyle \frac{\EZ{V_{m,n}^{X^{d_1}}}}{\EZ{V_{m,n}}}  \\
			& = &  
			\displaystyle\frac{\frac{1}{m} \sum_{i=1}^{m} s^2_n(X_i,\tX_i) + m_n(X_i)m_n(\tX_i)
			-
			 \frac{1}{m^2}\sum_{i,j=1}^{m} s^2_n(X_i,\tX_j) + m_n(X_i)m_n(\tX_j)
			}{
			\frac{1}{m} \sum_{i=1}^{m} s^2_n(X_i,X_i) + m_n(X_i)m_n(X_i)
			-
			 \frac{1}{m^2}\sum_{i,j=1}^{m} s^2_n(X_i,X_j) + m_n(X_i)m_n(X_j)
			} 
			\end{array}
			\end{equation}
			
			We note that the expression of $\tilde{\mathcal{S}}_{m,n}^{X^{d_1}} $ is different from the one obtained by  estimating ${V_m^{X^{d_1}}}/{V_m}$ in (\ref{Sobolestim})   by replacing $z(x)$ by the predictive mean  $m_n(x)$.  In $\tilde{\mathcal{S}}_{m,n}^{X^{d_1}}$ we take into account the kriging predictive covariance through the terms $s^2_n(X_i,\tX_j)$ and $s^2_n(X_i,X_j)$.

		\section{Kriging-based sensitivity analysis: a second approach}\label{chap6sec3}

		We have highlighted at the end of Subsection \ref{oakleymethod} that one of the main flaws  of the method presented by \cite{Oak04} is that it does not care about the exact definition of Sobol indices. We present in Subsection \ref{SobolEstimation2} another approach which  deals with this issue. Then, in Subsection    \ref{chap6sec30} we present an efficient method to compute  it.

			\subsection{Kriging-based Sobol index estimation}\label{SobolEstimation2}
		
			First of all, in  the previous section we have considered the variance of the main effects $V^{X^{d_1}}$ and  the total variance $V$ separately  in equation  (\ref{VX1V}). That is why the ratio of the expectations is considered as a sensitivity index in  equation (\ref{meanoak}). In fact, in a Sobol index  framework, we are interested in the ratio between $V^{X^{d_1}}$ and  $V$. Therefore, we suggest to deal directly with the following estimator (see equation (\ref{SobolestimZ})):
			\begin{equation}\label{point} 
			\mathcal{S}_{m,n}^{X^{d_1}} = \frac{V_{m,n}^{X^{d_1}}}{V_{m,n}}
			\end{equation}
			which corresponds to the ratio 	$V^{X^{d_1}}/V$  after substituting the code $z(x)$ by the Gaussian process $Z_n(x)$ and estimating the terms ${\VX{\EX{Z_n(X)|X^{d_1}}}}$ and ${\VX{Z_n(X)}}$ with a Monte-Carlo procedure as presented in \cite{sobol1993}. We note that we can naturally adapt the presented estimator with the ones  suggested by \cite{sobol2007estimating} and  \cite{Jan12}. 
			Nevertheless, we cannot obtain closed form expressions for the mean or the variance of this estimator. We thus have to numerically estimate them. 
			We present in Algorithm \ref{algoSobol1}  the suggested method to compute the distribution of  $\mathcal{S}_{m,n}^{X^{d_1}}$.

			\begin{algorithm}
			\caption{Evaluation of the distribution of $\mathcal{S}_{m,n}^{X^{d_1}}$.}
			\label{algoSobol1}
			\begin{algorithmic}[1]
			\State Build $Z_n(x)$ from the $n$ observations $\zn$ of $z(x)$ at points in $\D$ (see equation (\ref{predictiveZn})).
			\State  Generate two samples $(x_i)_{i=1,\dots,m}$ and $(\tx_i)_{i=1,\dots,m}$  of the  random vectors $(X_i)_{i=1,\dots,m}$ and $(\tX_i)_{i=1,\dots,m}$ with respect to the probability  measure $\mu$ (see Proposition \ref{varcovsob}).
			\State Set $N_Z$ the number of samples for $Z_n(x)$ and $B$ the number of bootstrap samples for evaluating the uncertainty due to Monte-Carlo  integrations.
			\For { $k=1,\dots,N_Z$}
				\State Sample   a realization $z_n(\xx)$ of $Z_n(\xx)$  with $\xx= \{(x_i)_{i=1,\dots,m}, (\tx_i)_{i=1,\dots,m} \} $
				\State Compute $\hat{\mathcal{S}}_{m,n,k,1}^{X^{d_1}}$ thanks to the equation (\ref{SobolestimZ}) from $z_n(\xx )$. 
				\For {l=2,\dots,B}
					\State Sample with replacements two  samples $\mathbf{u}$ and $\tilde{\mathbf{u}}$ from $\{(x_i)_{i=1,\dots,m}\} $ and $\{(\tx_i)_{i=1,\dots,m}\}$.
					\State Compute $\hat{\mathcal{S}}_{m,n,k,l}^{X^{d_1}}$ from $z_n(\xx^B)$ with $\xx^B = \{\mathbf{u}, \tilde{\mathbf{u}}\}$.
				\EndFor
			\EndFor

			\Return $\left(\hat{\mathcal{S}}_{m,n,k,l}^{X^{d_1}}\right)_{\substack{k=1,\dots,N_Z \\ l = 1,\dots,B}}$
			\end{algorithmic}
			\end{algorithm}
			
			The output $\left(\hat{\mathcal{S}}_{m,n,k,l}^{X^{d_1}}\right)_{\substack{k=1,\dots,N_Z \\ l = 1,\dots,B}}$ of Algorithm \ref{algoSobol1}  is a sample of size $N_Z \times B$   of  $\mathcal{S}_{m,n}^{X^{d_1}}$ defined on $(\Omega_X \times \Omega_Z, \sigma(\mathcal{F}_X\times \mathcal{F}_Z), \mathbb{P}_X \times \mathbb{P}_Z)$ (i.e. $\mathcal{S}_{m,n}^{X^{d_1}}$  takes both into account the uncertainty of the metamodel and the one of the Monte-Carlo integrations). Then, we can deduce the following estimate $\bar{\mathcal{S}}_{m,n}^{X^{d_1}}$ for  $\mathcal{S}_{m,n}^{X^{d_1}}$:
\begin{equation}\label{meanestimatorS}
\bar{\mathcal{S}}_{m,n}^{X^{d_1}} = \frac{1}{N_Z B} \sum_{\substack{k=1,\dots,N_Z \\ l = 1,\dots,B}}
\hat{\mathcal{S}}_{m,n,k,l}^{X^{d_1}}
\end{equation}
Furthermore, we can estimate the variance  of $\mathcal{S}_{m,n}^{X^{d_1}}$ with
\begin{equation}
\hat{\sigma}^2(\mathcal{S}_{m,n}^{X^{d_1}}) = \frac{1}{N_ZB-1} \sum_{\substack{k=1,\dots,N_Z \\ l = 1,\dots,B}}
\left( \hat{\mathcal{S}}_{m,n,k,l}^{X^{d_1}} - \bar{\mathcal{S}}_{m,n}^{X^{d_1}} \right)^2
\end{equation}

We note that the computational limitation of the algorithm is the sampling    of the Gaussian process  $Z_n(x)$   on $ \xx = \{(x_i)_{i=1,\dots,m}, (\tx_i)_{i=1,\dots,m} \} $. For that reason, we use a bootstrap procedure to evaluate  the uncertainty of the Monte-Carlo integrations instead of sampling   different realizations of the  random vectors $(X_i)_{i=1,\dots,m}$ and $(\tX_i)_{i=1,\dots,m}$.   Furthermore, the same  bootstrap samples are used for the $N_Z$ realizations of $Z_n(x)$.

Nevertheless, the number of Monte-Carlo particles $m$ is very large in general - it is often   around $m=5000d$ - and it thus can be an issue to compute realizations of $Z_n(x)$ on $\xx$. We present in the   Subsection  \ref{chap6sec30} an efficient method to deal with this point for any choice of $\mu(x)$, $\ff(x)$ and $r(x,\tx)$ and any index order.
The idea to carry out an estimation of (\ref{point}) from realizations of conditional Gaussian processes has already been suggested in \cite{gramacy2012categorical}. The main contribution  of this section is  the procedure  to balance the Monte-Carlo and the meta-model errors (see  Subsection \ref{optimalm}).

			\subsection{Determining the minimal number of Monte-Carlo particles $m$}\label{optimalm}
			We are interested here in quantifying  the uncertainty of the considered estimator $\mathcal{S}_{m,n}^{X^{d_1}}$ (\ref{point}). This estimator integrates two sources of uncertainty, the first one is related to the meta-model approximation and the second one is related to the Monte-Carlo integration. Therefore, we can decompose the variance of $\mathcal{S}_{m,n}^{X^{d_1}}$ as follows:
			\begin{displaymath}
			\mathrm{var}\left( \mathcal{S}_{m,n}^{X^{d_1}} \right) =  \mathrm{var}_Z\left(  \mathbb{E}_X \left[\mathcal{S}_{m,n}^{X^{d_1}} \big| Z_n(x)  \right] \right)
+
\mathrm{var}_X\left(  \mathbb{E}_Z \left[\mathcal{S}_{m,n}^{X^{d_1}} \big| (X_i, \tilde X_i)_{i=1,\dots,m}  \right] \right)
			\end{displaymath}
			where $\mathrm{var}_Z\left(  \mathbb{E}_X \left[\mathcal{S}_{m,n}^{X^{d_1}} \big| Z_n(x)  \right] \right)$ is the contribution of the meta-model on the variability of $\mathcal{S}_{m,n}^{X^{d_1}}$ and $\mathrm{var}_X\left(  \mathbb{E}_Z \left[\mathcal{S}_{m,n}^{X^{d_1}} \big| (X_i, \tilde X_i)_{i=1,\dots,m}  \right] \right)$ is the one of the Monte-Carlo integration.
			Furthermore, we have the following equalities:
			\begin{displaymath}
			\left\{ 
			\begin{array}{lll}
			\mathrm{var}_Z\left(  \mathbb{E}_X \left[\mathcal{S}_{m,n}^{X^{d_1}} \big| Z_n(x)  \right] \right) 
			& = & 
\mathbb{E}_X\left[  \mathrm{var}_Z \left(\mathcal{S}_{m,n}^{X^{d_1}} \big| (X_i, \tilde X_i)_{i=1,\dots,m}  \right)\right]
			\\
\mathrm{var}_X\left(  \mathbb{E}_Z \left[\mathcal{S}_{m,n}^{X^{d_1}} \big| (X_i, \tilde X_i)_{i=1,\dots,m}  \right] \right)
			 & = &
			\mathbb{E}_Z\left[   \mathrm{var}_X \left(\mathcal{S}_{m,n}^{X^{d_1}} \big| Z_n(x)  \right) \right] 
			\\
			\end{array}
			\right.
			\end{displaymath}
			 Therefore, from the sample $\left(\hat{\mathcal{S}}_{m,n,k,l}^{X^{d_1}}\right)_{\substack{k=1,\dots,N_Z \\ l = 1,\dots,B}}$  we can  estimate the part of  variance of the  estimator $\mathcal{S}_{m,n}^{X^{d_1}}$  related  to the meta-modelling  as follows:
			\begin{equation}\label{uncertaintymodel}
			\hat{\sigma}^2_{Z_n}(\mathcal{S}_{m,n}^{X^{d_1}}) = \frac{1}{B} \sum_{l=1}^B  \frac{1}{N_Z-1}  \sum_{k=1}^{N_Z}  \left(\hat{\mathcal{S}}_{m,n,k,l}^{X^{d_1}} - \bar{\hat{\mathcal{S}}}_{m,n,l}^{X^{d_1}} \right)^2
			\end{equation}
			where
			$
			\bar{\hat{\mathcal{S}}}_{m,n,l}^{X^{d_1}}  = \left( \sum_{i=1}^{N_Z} {\mathcal{S}}_{m,n,i,l}^{X^{d_1}}  \right)/ {N_Z } 
			$.
			Furthermore,  we can evaluate the  part of variance of $\mathcal{S}_{m,n}^{X^{d_1}}$ related  to the Monte-Carlo integrations as follows:
			\begin{equation}\label{uncertaintyMC}
			\hat{\sigma}^2_{MC}(\mathcal{S}_{m,n}^{X^{d_1}})  =\frac{1}{N_Z} \sum_{i=1}^{N_Z}   \frac{1}{B-1}  \sum_{i=1}^{B}  \left(\hat{\mathcal{S}}_{m,n,k,i}^{X^{d_1}} - \bar{\bar{\hat{\mathcal{S}}}}_{m,n,k}^{X^{d_1}} \right)^2
			\end{equation}
			where $\bar{\bar{\hat{\mathcal{S}}}}_{m,n,k}^{X^{d_1}}   = \left( \sum_{i=1}^{B} {\mathcal{S}}_{m,n,k,i}^{X^{d_1}}  \right)/ {B }$. 

			Therefore, we  have three different cases:
			\begin{enumerate}
			\item $\hat{\sigma}^2_{Z_n}(\mathcal{S}_{m,n}^{X^{d_1}})  \gg \hat{\sigma}^2_{MC}(\mathcal{S}_{m,n}^{X^{d_1}}) $: 
			the estimation error of $\mathcal{S}_{m,n}^{X^{d_1}}$ is essentially due to  the metamodel error.
			\item $\hat{\sigma}^2_{Z_n}(\mathcal{S}_{m,n}^{X^{d_1}})  \ll \hat{\sigma}^2_{MC}(\mathcal{S}_{m,n}^{X^{d_1}}) $: the estimation error of $\mathcal{S}_{m,n}^{X^{d_1}}$ is essentially due to the  Monte-Carlo error.
			\item $\hat{\sigma}^2_{Z_n}(\mathcal{S}_{m,n}^{X^{d_1}})  \approx \hat{\sigma}^2_{MC}(\mathcal{S}_{m,n}^{X^{d_1}}) $: the metamodel and the Monte-Carlo errors have the same contribution on the estimation error of $\mathcal{S}_{m,n}^{X^{d_1}}$.
			\end{enumerate}
			Considering that the number of observations $n$ is fixed, the minimal number of Monte-Carlo particles  $m$ is the one such that $\hat{\sigma}^2_{Z_n}(\mathcal{S}_{m,n}^{X^{d_1}})  \approx \hat{\sigma}^2_{MC}(\mathcal{S}_{m,n}^{X^{d_1}}) $. 
			We call it ``minimal'' since it is the one from which the Monte-Carlo error no  longer dominates. Therefore, it should be   the minimum number of required particles in   practical applications.
			In practice, to determine it,   we start with a small value of  $m$  and we increase it  while  the inequality $\hat{\sigma}^2_{Z_n}(\mathcal{S}_{m,n}^{X^{d_1}})  > \hat{\sigma}^2_{MC}(\mathcal{S}_{m,n}^{X^{d_1}}) $ is true. 
	
			\subsection{Sampling with respect to the kriging  predictive distribution on large data sets}\label{chap6sec30}

			We saw in the previous subsection in Algorithm \ref{algoSobol1}  that in a kriging framework, we can assess the distribution of the Sobol index estimators  from realizations of the  conditional Gaussian process  $Z_n(x)$ at points in $\xx$. Nevertheless, the size of the corresponding random  vector could be important since it equals twice the number of Monte-Carlo particles $m$. Therefore, computing such realizations could lead to numerical issues such as ill-conditioned matrix or huge computational cost, especially if we use a Cholesky decomposition. Indeed,  Cholesky decomposition complexity is $\mathcal{O}((2m)^3)$ and it often leads to ill-conditioned matrix since the predictive variance of $Z_n(x)$ is close to zero around the experimental design points.
 			
			 Let us introduce the following unconditioned Gaussian process:
			\begin{equation}
			\tilde{Z}(x) \sim \mathrm{GP}(0,\sigma^2r(x,\tx))
			\end{equation}
			We have the following proposition \cite{chiles1999geostatistics}:
			\begin{prop}[Sampling $Z_n(x)$ by kriging conditioning]\label{propsampZ}
			Let us consider the following Gaussian process:
			\begin{equation}
			\tilde{Z}_n(x) = m_n(x) - \tilde{m}_n(x) + \tilde{Z}(x)
			\end{equation}
			where $m_n(x)$ is the predictive mean of $Z_n(x)$ (\ref{predictiveZn}), 
			\begin{equation}
			 \tilde{m}_n(x)  = \ff'(x) \tilde{\bbeta} + \rr'(x)\RR^{-1} \left( \tilde{Z}(\D) - \F \tilde{\bbeta}\right)
			\end{equation}
			and $\tilde{\bbeta} = \left( \F'\RR^{-1}\F\right)^{-1}\F'\RR^{-1} \tilde{Z}(\D) $. Then, we have   
			\begin{displaymath}
			\tilde{Z}_n(x) \eqL Z_n(x)
			\end{displaymath}
			where $Z_n(x)$ has the distribution of the Gaussian process $Z(x)$ of mean $\ff'(x) \bbeta$ and covariance kernel $\sigma^2r(x,\tx)$ conditioned by $\zn$ at points in $\D$ (\ref{predictiveZn}). We note that we are in a Universal kriging case, i.e. we infer from the parameter $\bbeta$. In a simple kriging case, the proposition remains true by setting $\tilde{\bbeta} = 0$.
			\end{prop}

			The strength of Proposition \ref{propsampZ} is that it allows for sampling with respect to the distribution of $Z_n(x)$ by sampling an unconditioned Gaussian process $\tilde{Z}(x)$. The first consequence is that the conditioning of the covariance matrix is better since the variance of $\tilde{Z}(x)$ is not close to  zero around  points in $\D$. The second important consequence is that it allows for using efficient algorithms to compute realizations of $\tilde{Z}(x)$. For example, if $r(x,\tx)$ is a stationary kernel, one can use  the Bochner's theorem  (\cite{S99} p.29) and  the Fourier representation of $\tilde{Z}(x)$  to compute realizations of $\tilde{Z}(x)$ as presented in \cite{S99}.  Furthermore, another  efficient method is to use the Mercer's representation of $r(x,\tx)$ (see \cite{konig1986eigenvalue} and \cite{ferreira2009eigenvalues}) and the  Nystr\"om procedure  to approximate the  Karhunen-Loeve decomposition of $\tilde{Z}(x)$ as presented in \cite{R06} p.98. One of the main advantage of the Karhunen-Loeve decomposition of $Z(x)$ is that it allows for sequentially adding new points to $\xx$ without re-estimating the decomposition. Therefore, we can easily obtain the values of a given realization $z_n(x)$ of $Z_n(x)$ at new points not in $\xx$. This interesting property will allow us to efficiently estimate the number $m$ of Monte-Carlo particles such that the metamodel error and the Monte-Carlo estimation one are equivalent (see Subsection \ref{optimalm}).

		\section{Multi-fidelity co-kriging based sensitivity analysis}\label{chap6sec4}

		Now let us suppose that we have $s$ levels of code $(z_t(x))_{t=1,\dots,s}$ from the less accurate one $z_1(x)$ to the most accurate one $z_s(x)$ and  that we want to perform a Global sensitivity analysis for  $z_s(x)$. 
We consider that, conditioning on the model parameters, $(z_t(x))_{t={1,\dots,s}}$ are realizations   of Gaussian processes $(Z_t(x))_{t={1,\dots,s}}$. Furthermore, we consider the following multi-fidelity model $t=2,\dots,s$:
		\begin{equation}\label{mufimodel}
		\left\{
		\begin{array}{l}
		Z_t(x) = \rho_{t-1} Z_{t-1}^*(x) + \delta_t(x) \\
		Z_{t-1}^*(x) \perp  \delta_t(x) \\
		Z_{t-1}^*(x)  \sim [Z_{t-1}(x) |\ZZ^{(t-1)} = \zz^{(t-1)}, \bbeta, \rrho, \ssigma^2  ]
		\end{array}
		\right.
		\end{equation}
		where  $\bbeta = (\bbeta_t)_{t=1,\dots,s}$, $\rrho = (\rho_{t-1})_{t=2,\dots,s}$,  $\ssigma^2 = (\sigma_t^2)_{t=1,\dots,s}$, $\ZZ^{(l-1)}  = (Z_1(\D_1),\dots, Z_{t-1}(\D_{t-1}))$, $\zz^{(t-1)}  = (z_1(\D_1),\dots, z_{t-1}(\D_{t-1}))$ and $(\D_t)_{t=1,\dots,s}$ are the experimental design sets  at level $t$ with $n_t$ points and such that $\D_s \subseteq \D_{s-1} \subseteq \dots \subseteq \D_1$. Further, conditioning on $\bbeta_t$ and $\sigma^2_t$, $\delta_t(x) $ is a Gaussian process of mean $\ff'_t(x) \bbeta_t$ and covariance $\sigma_t^2r_t(x,\tx)$ and we use the convention $Z_1(x) = \delta_1(x)$.
		This model is analogous to the one presented in \cite{KO00} except that $Z_{t-1}^*(x)$ has a conditional distribution. 

		We propose a Bayesian formulation of the model which allows to consider non-informative prior distributions for the  the regression  parameters $\bbeta = (\bbeta_t)_{t=1,\dots,s}$ and the adjustment parameters $\rrho = (\rho_{t-1})_{t=2,\dots,s}$.  This leads to the following predictive distribution which  integrates the posterior distributions of the parameters  $\bbeta = (\bbeta_t)_{t=1,\dots,s}$ and $\rrho = (\rho_{t-1})_{t=2,\dots,s}$.
		\begin{equation}\label{predmufisobol}
		[Z_s(x)|\ZZ^{(s)} = \zz^{(s)}, \ssigma^2]
		\end{equation}
		The predictive distribution (\ref{predmufisobol})  is not    Gaussian.  Nevertheless, we can have closed form expressions for  its  mean $\mu_{n_s}^s(x)$ and covariance $k_{n_s}^s(x,\tx)$:
		\begin{equation}\label{mufisobolmean}
		\mu_{n_s}^s(x) = \hat{\rho}_{s-1}\mu_{n_{s-1}}^{s-1}(x) + \mu_{\delta_s}(x) 
		\end{equation}
		and:
		\begin{equation}\label{mufisobolvar}
		k_{n_s}^s(x,\tx) = \widehat{\rho_{s-1}^2}k_{n_{s-1}}^{s-1}(x,\tx) +k_{\delta_s}(x,\tx)
		\end{equation}
		where for 
  $t=1,\dots,s$,  $\begin{pmatrix} \hat{\rho}_{t-1} \\  \hat{\bbeta}_t \end{pmatrix} = ( \HHH_t' \RR_t^{-1}\HHH_t)^{-1}\HHH_t' \RR_t^{-1}\zz_t  $, $\HHH_t =[z_{t-1}(\D_t) \quad \F_t] $, $\F_t = \ff'_t(\D_t)$, $\hat \rho_0 = 0$,  $\HHH_1 = \F_1$, $\widehat{\rho_{t-1}^2} = \hat{\rho}_{t-1}^2 +  \left[ ( \HHH_t' \RR_t^{-1}\HHH_t)^{-1}\right]_{[1,1]} $,  $\widehat{\rho_{0}^2} = 0$, 
		\begin{equation}\label{mudeltat}
		\mu_{\delta_t}(x) =  \ff'_t(x)\hat{\bbeta}_t +     \rr_t'(x)\RR_t^{-1}(\zz_t-\F_t\hat{\bbeta}_t - \hat{\rho}_{t-1}z_{t-1}(\D^{t}))
		\end{equation}
		and
		\begin{equation}\label{kdeltat}
		k_{\delta_t}(x,\tx) = \sigma_t^2 \left( r_t(x,\tx)  -  \begin{pmatrix} \hh_t'(x) &  \rr_t' (x)  \end{pmatrix}   \begin{pmatrix} 0 & \HHH_t' \\ \HHH_t   & \RR_t  \end{pmatrix} ^{-1} \begin{pmatrix} \hh_t(\tx) \\ \rr_t(\tx) \end{pmatrix} \right)
		\end{equation}
		with  $\hh_t'(x)  = [\mu_{n_{t-1}}^{t-1}(x) \quad \ff'_t(x)]$ and $\hh_1'(x)  = \ff'_1(x)$. We note that, in the mean of the predictive distribution, the regression and adjustment  parameters have  been replaced by their  posterior means. Furthermore, the predictive variance integrates the uncertainty due to the regression and adjustment parameters.
		
 We note that for each $t=1,\dots,s$,  the variance parameter  $\sigma_t^2$ is estimated with a restricted maximum likelihood method. Thus, its  estimation is given by $\hat{\sigma}_t^2 = (\zz_t - \HHH_t\hat{\bbeta}_t)'\RR_t^{-1}(\zz_t - \HHH_t\hat{\bbeta}_t)/(n_t - p_t - 1)$ where $p_t$ is the size of $\bbeta_t$.

		We present in Subsection \ref{chap6sec41} the extension   in a multi-fidelity framework of the     first kriging-based Sobol index estimations presented in \cite{Oak04}. Then, we present in Subsection \ref{chap6sec40} the extension of our approach to perform co-kriging-based multi-fidelity sensitivity analysis.

			\subsection{Extension   of the first approach  for multi-fidelity co-kriging models}\label{chap6sec41}
			Let us denote by $\tilde{\mathcal{S}}_{m,s}^{X^{d_1}}$ the estimation of
  $  V^{X^{d_1}}  /  V  $ when we substitute $z_s(x)$ by $Z_{n,s}(x) \sim [Z_s(x)|\ZZ^{(s)} = \zz^{(s)}, \sigma^2]$ and when we use the Sobol procedure to perform Monte-Carlo estimations (see \cite{sobol1993} and Subsection \ref{MCsobolestim}). Then, the estimator suggested in \cite{Oak04} and \cite{Mar09} becomes in a multi-fidelity framework:
			\small
			\begin{eqnarray*} 
			\tilde{\mathcal{S}}_{m,s}^{X^{d_1}} 
			& = &  \frac{\frac{1}{m} \sum_{i=1}^{m} k_{n_s}^s(X_i,\tX_i) + \mu_{n_s}^s(X_i)\mu_{n_s}^s(\tX_i)
			-
			 \frac{1}{m^2}\sum_{i,j=1}^{m} k_{n_s}^s(X_i,\tX_j) + \mu_{n_s}^s(X_i)\mu_{n_s}^s(\tX_j)
			}{
			\frac{1}{m} \sum_{i=1}^{m} k_{n_s}^s(X_i,X_i) + \mu_{n_s}^s(X_i)\mu_{n_s}^s(X_i)
			-
			 \frac{1}{m^2}\sum_{i,j=1}^{m} k_{n_s}^s(X_i,X_j) + \mu_{n_s}^s(X_i)\mu_{n_s}^s(X_j)
			}  \\
			& = & \frac{U}{D} 
			\end{eqnarray*}
			\normalsize
			where 
			\begin{eqnarray*}
			U & =   & 
			\frac{1}{m} \sum_{i=1}^{m} \left( \sum_{t=1}^s \left( \prod_{j=t}^{s-1}\widehat{\rho^2_j }\right) k_{\delta_t} (X_i,\tX_i) + 
			 \sum_{t,\tildt=1}^s \left( \prod_{j=t}^{s-1}\hat{\rho}_j \right) \left( \prod_{j=\tildt}^{s-1}\hat{\rho}_j \right)
			\mu_{\delta_t} (X_i)   \mu_{\delta_\tildt} (\tX_i)
			\right)
			\\
			 & & -
			 \frac{1}{m^2}\sum_{i,j=1}^{m}  \left( \sum_{t=1}^s \left( \prod_{j=t}^{s-1}\widehat{\rho^2_j }\right) k_{\delta_t} (X_i,\tX_j) + 
			 \sum_{t,\tildt=1}^s \left( \prod_{j=t}^{s-1}\hat{\rho}_j \right) \left( \prod_{j=\tildt}^{s-1}\hat{\rho}_j \right)
			\mu_{\delta_t} (X_i)   \mu_{\delta_\tildt} (\tX_j)
			\right)
			\end{eqnarray*}
			\begin{eqnarray*}
			D &  = &   
			\frac{1}{m} \sum_{i=1}^{m} \left( \sum_{t=1}^s \left( \prod_{j=t}^{s-1}\widehat{\rho^2_j} \right) k_{\delta_t} (X_i,X_i) + 
			 \sum_{t,\tildt=1}^s \left( \prod_{j=t}^{s-1}\hat{\rho}_j \right) \left( \prod_{j=\tildt}^{s-1}\hat{\rho}_j \right)
			\mu_{\delta_t} (X_i)   \mu_{\delta_\tildt} (X_i)
			\right)
			\\
			& & -
			 \frac{1}{m^2}\sum_{i,j=1}^{m}  \left( \sum_{t=1}^s \left( \prod_{j=t}^{s-1}\widehat{\rho^2_j }\right) k_{\delta_t} (X_i,X_j) + 
			 \sum_{t,\tildt=1}^s \left( \prod_{j=t}^{s-1}\hat{\rho}_j \right) \left( \prod_{j=\tildt}^{s-1}\hat{\rho}_j \right)
			\mu_{\delta_t} (X_i)   \mu_{\delta_\tildt} (X_j)
			\right)
			\end{eqnarray*}
			and with the conventions $\hat{\rho}_0 = 0$, $\prod_{i=s}^{s-1} \hat{\rho}_i = 1$, $\widehat{\rho^2_0} = 0$, $\prod_{i=s}^{s-1}\widehat{\rho^2_i} = 1$,  $\mu_{\delta_1} (x)  = \mu_{n_1}^1(x) $ and $k_{\delta_1} (x,\tx)  = k_{n_1}^1(x,\tx)$.

			We note that $\tilde{\mathcal{S}}_{m,s}^{X^{d_1}}$ is the analogous of $\tilde{\mathcal{S}}_{m,n}^{X^{d_1}}$ presented in Subsection \ref{MCoak}.
			Furthermore,  the developed expression  of $\tilde{\mathcal{S}}_{m,s}^{X^{d_1}}$ allows for identifying the contribution of each code level $t$ to the sensitivity index  and the one of the covariance between the bias and the code at level $t$. We note that the covariance here is with respect to the distribution of the input parameters $X$.
			Nevertheless, as pointed out in previous sections, this estimator is based  on a ratio  of expectations and  thus does not correspond to the true Sobol indices. 

			\subsection{Extension   of the second approach  for multi-fidelity co-kriging models}\label{chap6sec40}

			We present here  the extension of the approach presented in Section \ref{chap6sec3} to the multi-fidelity co-kriging model. Therefore, we aim to sample with respect to the distribution of 
			\begin{equation} \label{estimatormufi}
			\mathcal{S}_{m,s}^{X^{d_1}} =  \frac{\frac{1}{m} \sum_{i=1}^{m} Z_{n,s}(X_i)Z_{n,s}(\tX_i)
			-
			 \frac{1}{m}\sum_{i=1}^{m} Z_{n,s}(X_i)  \frac{1}{m}\sum_{i=1}^{m}Z_{n,s}(\tX_i) 
			}{
			\frac{1}{m} \sum_{i=1}^{m} Z_{n,s}(X_i)^2
			-
			\left( \frac{1}{m}\sum_{i=1}^{m} Z_{n,s}(X_i)\right)^2
			} 
			\end{equation}
			which is the analog of $\mathcal{S}_{m,n}^{X^{d_1}} $  (\ref{SobolestimZ}) in an univariate case when we substitute $z(x)$ with $Z_{n,s}(x) \sim [Z_s(x)|\ZZ^{(s)} = \zz^{(s)}, \sigma^2]$.  In fact, we can directly use Algorithm \ref{algoSobol1} by  sampling   realizations of  $Z_{n,s}(x) $ instead of $Z_n(x)$. Moreover,  the procedure presented in Subsection \ref{optimalm} to determine the optimal number of Monte-Carlo particles $m$ is straightforward.
	
			However, the distribution of $Z_{n,s}(x)$ is not Gaussian and thus the method  presented in Subsection \ref{chap6sec30} cannot be used directly.
			In order to handle this problem, we consider the conditional distribution   $ [Z_s(x)|\ZZ^{(s)} = \zz^{(s)}, \ssigma^2, \rrho, \bbeta]$, with $\ssigma^2 = (\sigma_t^2)_{t=1,\dots,s}$,  $\bbeta = (\bbeta_t)_{t=2,\dots,s}$ and $\rrho = (\rho_{t-1})_{t=2,\dots,s}$  which  is Gaussian (note that we infer from $\bbeta_1$).  
			It corresponds to the distribution (\ref{predmufisobol}) conditioning by $\bbeta$ and $\rrho$.
			Furthermore, the Bayesian estimation of $(\rho_{t-1}, \bbeta_t)$ gives us   for all $t=2,\dots,s$:
			\begin{equation}\label{samplebetarho}
			 \begin{pmatrix}  {\rho}_{t-1} \\   {\bbeta}_t \end{pmatrix} \sim \N{ ( \HHH_t' \RR_t^{-1}\HHH_t)^{-1}\HHH_t' \RR_t^{-1}\zz_t  }{\sigma_t^2( \HHH_t' \RR_t^{-1}\HHH_t)^{-1}}
			\end{equation}
			From the recursive formulation given in (\ref{mufimodel}),  we can  define  the following Gaussian process having  the desired distribution    $ [Z_s(x)|\ZZ^{(s)} = \zz^{(s)}, \ssigma^2, \rrho, \bbeta]$:
			\begin{equation}
			Z_{n,s,\rho,\bbeta}(x)  = \left(\prod_{j=1}^{s-1}\rho_{j} \right)  Z_{n,1}(x) + \sum_{t=2}^{s-1} \left(\prod_{j=t}^{s-1}\rho_{j} \right)  \delta_{t,\rho_{t-1},\bbeta_t}(x)+\delta_{s,\rho_{s-1},\bbeta_s}(x)
			\end{equation}
			 where (see equations (\ref{mudeltat}) and (\ref{kdeltat})):
			\begin{equation}\label{Zn1}
			Z_{n,1}(x)  \sim \mathrm{GP}\left(\mu_{\delta_1}(x),  k_{\delta_1}(x,\tx) \right)
			\end{equation}
			and for $t=2,\dots,s$:
			\begin{equation}\label{deltatrhobeta}
			 \delta_{t,\rho_{t-1},\bbeta_t}  (x)  \sim \mathrm{GP}\left(\mu_{t,\rho_{t-1},\bbeta_t}(x),  k_{t,\rho_{t-1},\bbeta_t}(x,\tx) \right)
			\end{equation}
			with 
			$
			\mu_{t,\rho_{t-1},\bbeta_t}(x) = \rr_t'(x)\RR_t^{-1}(\zz_t-\F_t{\bbeta}_t - {\rho}_{t-1}z_{t-1}(\D^{t}))
			$, $\left(( \delta_{t,\rho_{t-1},\bbeta_t}  (x) )_{t=2,\dots,s}, Z_{n,1}(x)  \right)$ independent 
			and
			$$
			k_{t,\rho_{t-1},\bbeta_t}(x,\tx) = \sigma_t^2 \left( r_t(x,\tx)  -     \rr_t' (x)  \RR_t ^{-1}  \rr_t(\tx)  \right).
			$$
			Therefore, we can deduce the following algorithm to compute a  realization $z_{n,s}(x)$  of 
$Z_{n,s}(x) \sim [Z_s(x)|\ZZ^{(s)} = \zz^{(s)}, \ssigma^2]$.
			\begin{algorithm}
			\caption{Sampling with respect to the predictive  distribution $[Z_s(x)|\ZZ^{(s)} = \zz^{(s)}, \sigma^2]$.}
			\label{algoSobol2}
			\begin{algorithmic}[1]
			\State Generate a sample $z_{n,1}(x)$  with respect to  (\ref{Zn1}) thanks to the method presented in Proposition \ref{propsampZ} in the   universal kriging case.
			\State Set $z_{n,s}(x) = z_{n,1}(x)$.
				\For {t=2,\dots,s}
					\State Generate a sample $ \begin{pmatrix}  {\rho}_{t-1}^* \\   {\bbeta}_t^*  \end{pmatrix} $  with respect to (\ref{samplebetarho}).
					\State Conditionally to $ \begin{pmatrix}  {\rho}_{t-1}^* \\   {\bbeta}_t^*  \end{pmatrix} $,  generate a sample  $\delta_{t,\rho_{t-1}^*,\bbeta_t^*}^*(x) $  with respect to (\ref{deltatrhobeta})  thanks to the method presented in Proposition \ref{propsampZ} in the simple kriging case.
					\State Set $z_{n,s}(x) = {\rho}_{t-1}^* z_{n,s} (x) +  \mu_{t,\rho_{t-1}^*,\bbeta_t^*}^*(x )$.	
				\EndFor

			\Return $z_{n,s}(x)$.
			\end{algorithmic}
			\end{algorithm}

			Algorithm \ref{algoSobol2}  provides an efficient tool to sample with respect to the distribution $[Z_s(x)|\ZZ^{(s)} = \zz^{(s)}, \sigma^2]$.  Then, from each sample  we can estimate the Sobol indices with a Monte-Carlo  procedure. Naturally, we can easily use a bootstrap procedure to take into account the uncertainty due to the Monte-Carlo  scheme.  Furthermore, we see in Algorithm \ref{algoSobol2} that once a  sample of $[Z_s(x)|\ZZ^{(s)} = \zz^{(s)}, \sigma^2]$ is available, a sample for each distribution $[Z_t(x)|\ZZ^{(t)} = \zz^{(t)}, \sigma^2]$, $t=1,\dots,s-1$ is also available. Therefore, we can directly quantify the difference between the Sobol indices at a level $t$ and the ones at another level $\tildt$.

\section{Numerical illustrations on an academic example}\label{sec:tests}

We illustrate here the  kriging-based sensitivity analysis suggested in Section \ref{chap6sec3}. We remind that the aim of this approach is to perform a sensitivity index taking into account both the uncertainty related to the surrogate modeling and the one related to the Monte-Carlo integrations. 
Let us  consider the Ishigami function:
\begin{displaymath}
z(x_1, x_2, x_3) =  \mathrm{sin}(x_1) + 7 \mathrm{sin(x_2)}^2+0.1x_3^4\mathrm{sin}(x_1),
\end{displaymath}
where $\mu_i$  is uniform on $[-\pi, \pi]$, $i=1,2,3$. We are interested in  the first order sensitivity indices theoretically given by
\begin{displaymath}
(S_1, S_2, S_3) = (0.314, 0.442, 0).
\end{displaymath}

This section is organized as follows. First, in Subsection \ref{compsob} we compare the Sobol index estimator  $\hat{\mathcal{S}}_{m,n}^{X^{d_1}}$  (\ref{meanoakestim})  proposed  by \cite{Oak04}, the suggested one  given by the mean of $\mathcal{S}_{m,n}^{X^{d_1}} $  (\ref{point}) and the usual one which consists in substituting $z(x)$ by the predictive mean $m_n(x)$ (\ref{predictiveZn}) in (\ref{Janonestim}). Then, in sections \ref{nincrease}, \ref{optimalmishi} and \ref{coverageMC} we deal with the approach presented in Section \ref{chap6sec3}. In particular, we show that this approach is relevant to perform an uncertainty quantification taking into account both  the uncertainty of the meta-modeling and the one of the Monte-Carlo integrations. We note that     the construction of the surrogate models   used in sections  \ref{nincrease}, \ref{optimalmishi} and \ref{coverageMC} is presented in Section \ref{krigbuilding}.

\subsection{Comparison between the different methods}\label{compsob}

\begin{figure}[!ht]
\begin{center}
\includegraphics[width = 7 cm, height = 7cm]{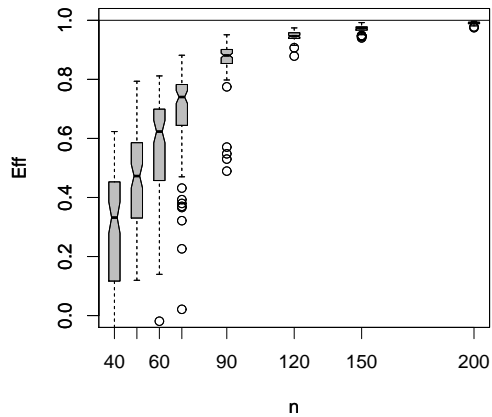}
\caption{Convergence of the model efficiency when the number $n$ of observations increases. 100  LHS  are randomly sampled for each number of observations $n$.}
\label{Q2effcompson}
\end{center}
\end{figure}

The aim of this subsection is to  perform a numerical comparison between   $\tilde{\mathcal{S}}_{m,n}^{X^{d_1}}$  (\ref{meanoakestim}), the  empirical mean of $\mathcal{S}_{m,n}^{X^{d_1}} $ given in Equation (\ref{meanestimatorS})  and the following estimator (see (\ref{Janonestim})):
			\begin{equation} \label{janonGPmean}
			\check{\mathcal{S}}_{m,n}^{X^{d_1}}  = \frac{\frac{1}{m} \sum_{i=1}^{m} m_n(X_i)m_n(\tX_i)
			-
			\left(  \frac{1}{2m}\sum_{i=1}^{m} m_n(X_i) + m_n(\tX_i)   \right)^2  
			}{
			\frac{1}{m} \sum_{i=1}^{m} m_n(X_i)^2
			-
			\left(  \frac{1}{2m}\sum_{i=1}^{m} m_n(X_i) + m_n(\tX_i)   \right)^2  
			} .
			\end{equation}
We note that the empirical mean $\bar{\mathcal{S}}_{m,n}^{X^{d_1}} $  of $\mathcal{S}_{m,n}^{X^{d_1}} $ is evaluated thanks to Algorithm  \ref{algoSobol1}, with $N_Z = 500$ and $B = 1$:
			\begin{equation*}
			\bar{\mathcal{S}}_{m,n}^{X^{d_1}}  = \frac{1}{N_Z }\sum_{ k=1,\dots,N_Z  }  \hat{\mathcal{S}}_{m,n,k,1}^{X^{d_1}},
			\end{equation*}
and for $\tilde{\mathcal{S}}_{m,n}^{X^{d_1}}$  and $\mathcal{S}_{m,n}^{X^{d_1}} $ we use the Monte-Carlo estimator (\ref{Janonestim}) suggested in \cite{Jan12} (it is the one used in (\ref{janonGPmean}).
Then for the comparison, different sizes of the learning sample are considered ($n=40,$  $50,$ $60,$ $70,$ $90,$ $120,$ $150,$ $200$ observations) and we  randomly build 100 Latin Hypercube Samples (LHS) for each size of the learning sample. From these experimental design sets, we build kriging models with a constant trend $\beta$ and a  tensorised  $5/2$-Mat\'ern kernel. Furthermore, the characteristic length scales   $(\theta_i)_{i=1,2,3}$ are estimated with a maximum likelihood procedure for each design set. The Nash-Sutcliffe model efficiency coefficient (sometimes called the predictivity coefficient $Q^2$), 
\begin{displaymath}
\mathit{Eff}_n = 1 - \frac{\sum_{x\in T} (m_n(x)- z(x))^2}{\sum_{x\in T} (m_n(x) - \bar{z}(x))^2},  \quad  \bar{z}(x) = \frac{1}{\#T} \sum_{x\in T}z(x),
\end{displaymath}
 of the different kriging models  are  evaluated on a test set $T$ composed of 1,000 points uniformly spread on the input parameter space $[-\pi, \pi]^3$. 
The values of $\mathit{Eff}_n$  are  presented in Figure \ref{Q2effcompson}.
The closer $\mathit{Eff}$ is to 1, the more accurate is the model $m_n(x)$.

\begin{figure}[!ht]
        \centering
        \begin{subfigure}[b]{0.45\textwidth}
                \centering
                \includegraphics[width = 6.5 cm, height = 5.5 cm]{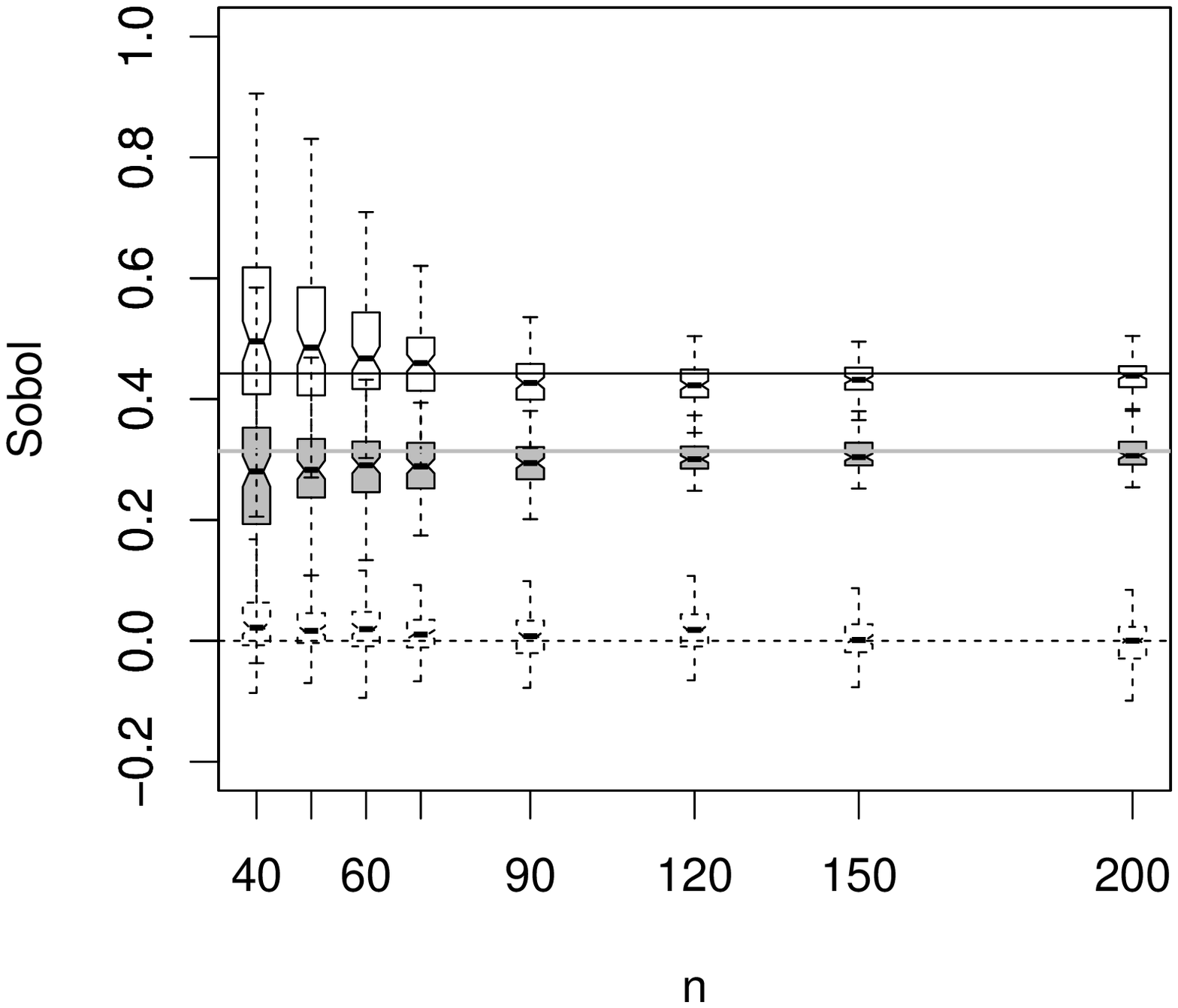}
	\captionsetup{format=hang,justification=centering}
                \caption[hang,center]{}
                \label{a:sobolcomparison}
        \end{subfigure}%
	~
        \begin{subfigure}[b]{0.45\textwidth}
                \centering
                \includegraphics[width = 6.5 cm, height = 5.5 cm]{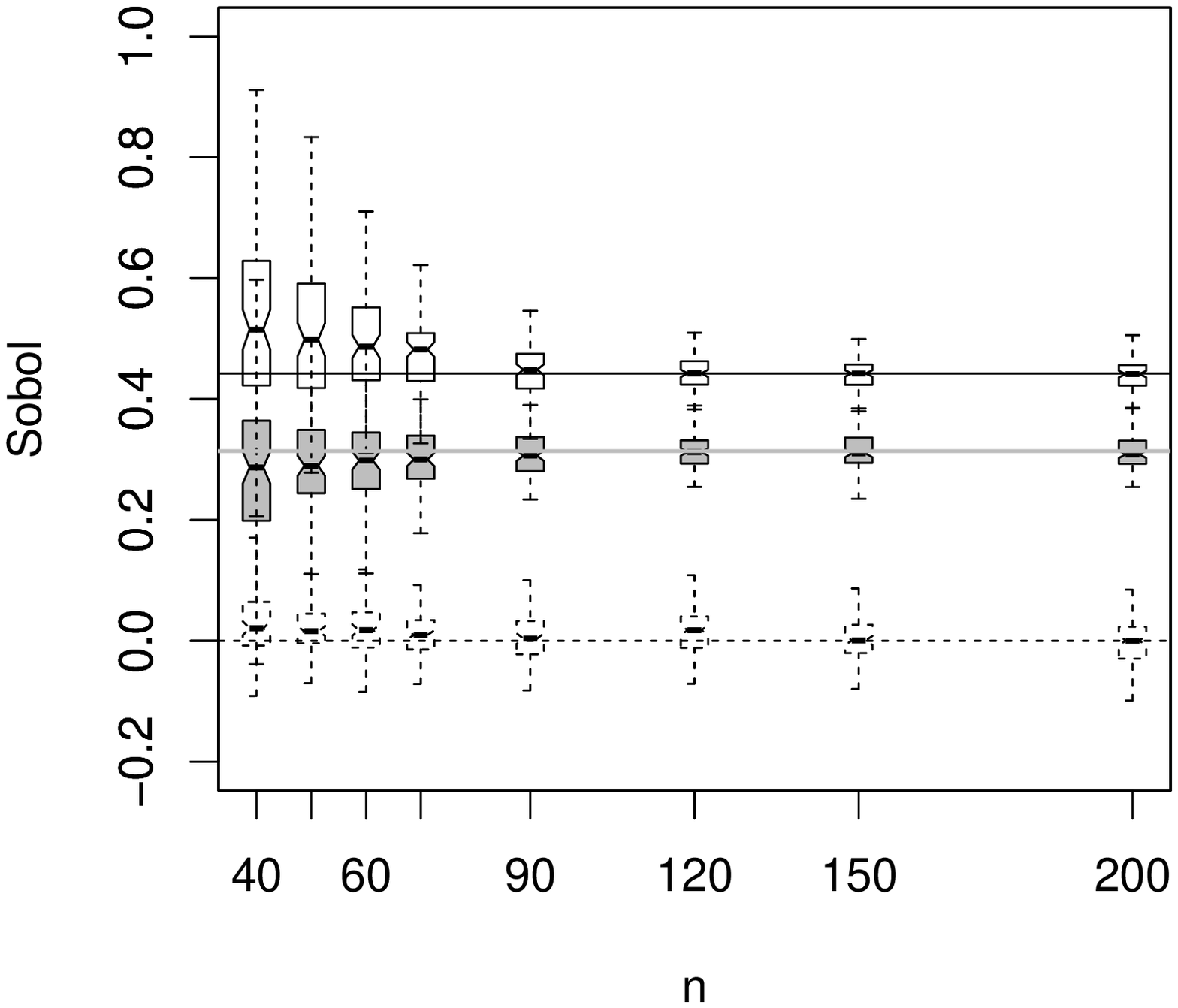}
	\captionsetup{format=hang,justification=centering}
                \caption[hang,center]{}
                \label{b:sobolcomparison}
        \end{subfigure}
 
        \begin{subfigure}[b]{0.45\textwidth}
                \centering               
                \includegraphics[width = 6.5 cm, height = 5.5 cm]{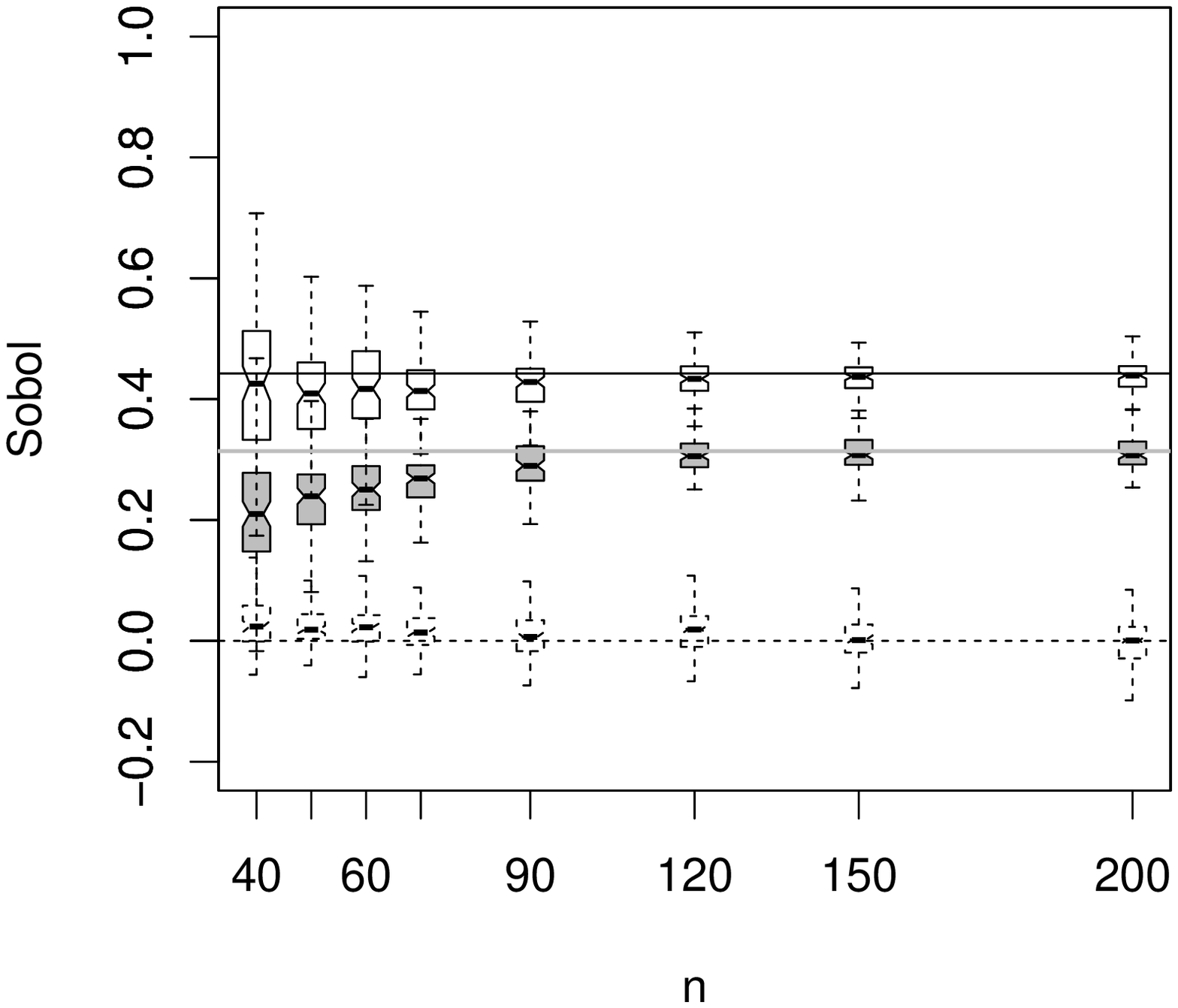}
	\captionsetup{format=hang,justification=centering}
                \caption[hang,center]{}
                \label{c:sobolcomparison}
        \end{subfigure}
      \caption{Comparison between three Sobol index estimators. The comparison are performed from 100 random LHS for each number of observations $n$. Figure (\subref{a:sobolcomparison}) corresponds to the suggested Sobol estimator $\bar{\mathcal{S}}_{m,n}^{X^{d_1}}$  (see Section \ref{chap6sec3}), Figure  (\subref{b:sobolcomparison})  corresponds to the usual (metamodel predictor only) estimator $\check{\mathcal{S}}_{m,n}^{X^{d_1}} $  (see Equation (\ref{janonGPmean})) and Figure  (\subref{c:sobolcomparison}) corresponds to the estimator $\tilde{\mathcal{S}}_{m,n}^{X^{d_1}}$   suggested in \cite{Oak04}. The horizontal lines represent the true values of the Sobol indices (solid gray line: $S_1$; solid black line: $S_2$ and dashed black line: $S_3$)}\label{sobolcomparison}
\end{figure}

Figure \ref{sobolcomparison} illustrates the Sobol index estimates obtained with the three methods.  We see in Figure \ref{sobolcomparison} that the suggested estimator  $\bar{\mathcal{S}}_{m,n}^{X^{d_1}}$ performs as well as the usual estimator $\check{\mathcal{S}}_{m,n}^{X^{d_1}} $ (\ref{janonGPmean}). In fact, as we will see in  the next subsections, the strength of the suggested estimator is  to  provide  more relevant uncertainty quantification. Finally, we see in Figure \ref{c:sobolcomparison} that the estimator $\tilde{\mathcal{S}}_{m,n}^{X^{d_1}}$  (\ref{meanoakestim}) suggested in \cite{Oak04} seems to systematically underestimate the true value of the Sobol index for non-negligible index and when the model efficiency is low.

\subsection{Model building and Monte-Carlo based estimator}\label{krigbuilding}

For the numerical illustrations  in sections  \ref{nincrease} and \ref{optimalmishi}, we use different kriging models built from different experimental design sets (optimized-LHS with respect to the centered $L_2$-discrepancy  criterion, \cite{damblin2013}) of   size $n=30,\dots,200$. Furthermore, for all kriging models, we consider a constant trend $\beta$ and a tensorised  $5/2$-Mat\'ern kernel (see \cite{R06}).

\begin{figure}[!ht]
\begin{center}
\includegraphics[width = 6 cm]{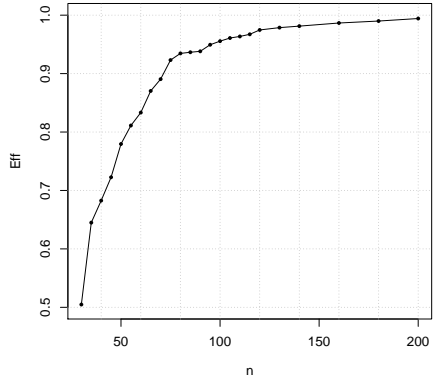}
\caption{Convergence of the model efficiency when the number $n$ of observations increases. For each number of observations $n$, the experimental design set is an optimized-LHS with respect to the centered $L_2$-discrepancy.}
\label{figQ2}
\end{center}
\end{figure}

The   characteristic length scales  $(\theta_i)_{i=1,2,3}$  are estimated for each experimental design set  by maximizing the marginal likelihood. Furthermore, the variance parameter $\sigma^2$ and the trend  parameter $\beta$ are   estimated  with   a maximum likelihood method for each experimental design set too.
Then for each $n$, the  Nash-Sutcliffe model efficiency  is evaluated on a test set composed of 10,000 points uniformly spread on the input parameter space $[-\pi, \pi]^3$. Figure  \ref{figQ2} illustrates the  estimated values of $\mathit{Eff}_n$  with respect to the number of observations $n$.

Then, for  estimating the Sobol indices, we use the Monte-Carlo based estimator given by (\ref{Janonestim}). It has the strength to be asymptotically efficient for  the first order indices (see \cite{Jan12}). 

\subsection{Sensitivity index estimates when $n$ increases}\label{nincrease}

Let us  consider a fixed number of Monte-Carlo particles $m=10,000$. The aim of this subsection is to quantify the part of the  index estimator uncertainty  related to the Monte-Carlo integrations and  
 the one related to the surrogate modeling.

To perform such analysis we use the procedure presented in Algorithm  \ref{algoSobol1}  with  $B = 300$ bootstrap samples   and $N_Z = 500$  realizations of $Z_n(x)$ (\ref{predictiveZn}). It results for each $i=1,2,3$ a sample $\left( \hat{S}^i_{m,n,k,l}\right)$, $k=1,\dots,N_Z$, $l=1,\dots,B$,  with respect to the distribution  of the estimator obtained by substituting $z(x)$ with $Z_n(x)$ in (\ref{Janonestim}).

\begin{figure}[!ht]
        \centering
        \begin{subfigure}[b]{0.45\textwidth}
                \centering
                \includegraphics[width = 6 cm, height = 5cm]{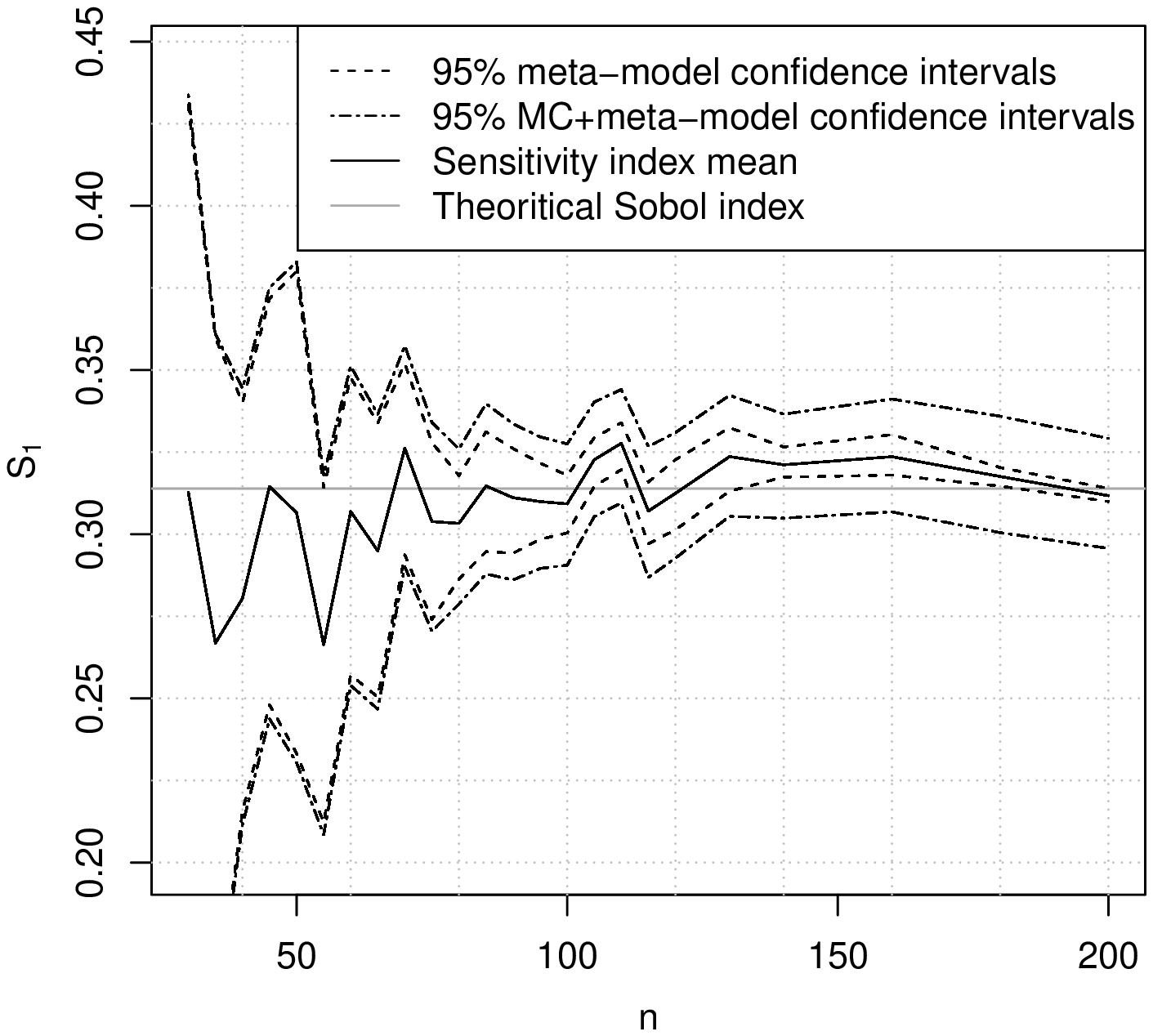}
	\captionsetup{format=hang,justification=centering}
                \caption[hang,center]{}
                \label{a:figS}
        \end{subfigure}%
	~
        \begin{subfigure}[b]{0.45\textwidth}
                \centering
                \includegraphics[width = 6 cm, height = 5cm]{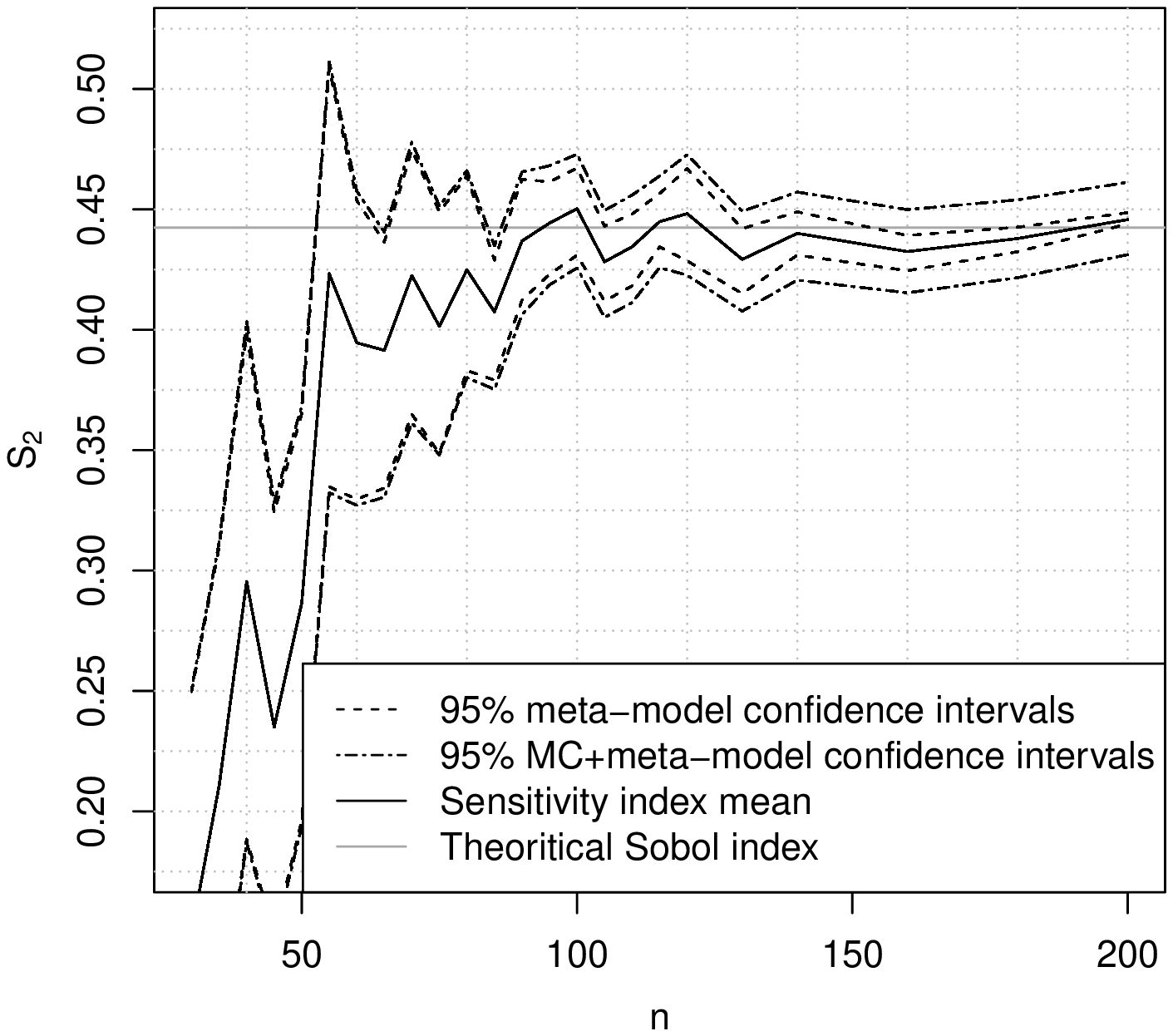}
	\captionsetup{format=hang,justification=centering}
                \caption[hang,center]{}
                \label{b:figS}
        \end{subfigure}
 
        \begin{subfigure}[b]{0.45\textwidth}
                \centering               
                \includegraphics[width = 6 cm, height = 5cm]{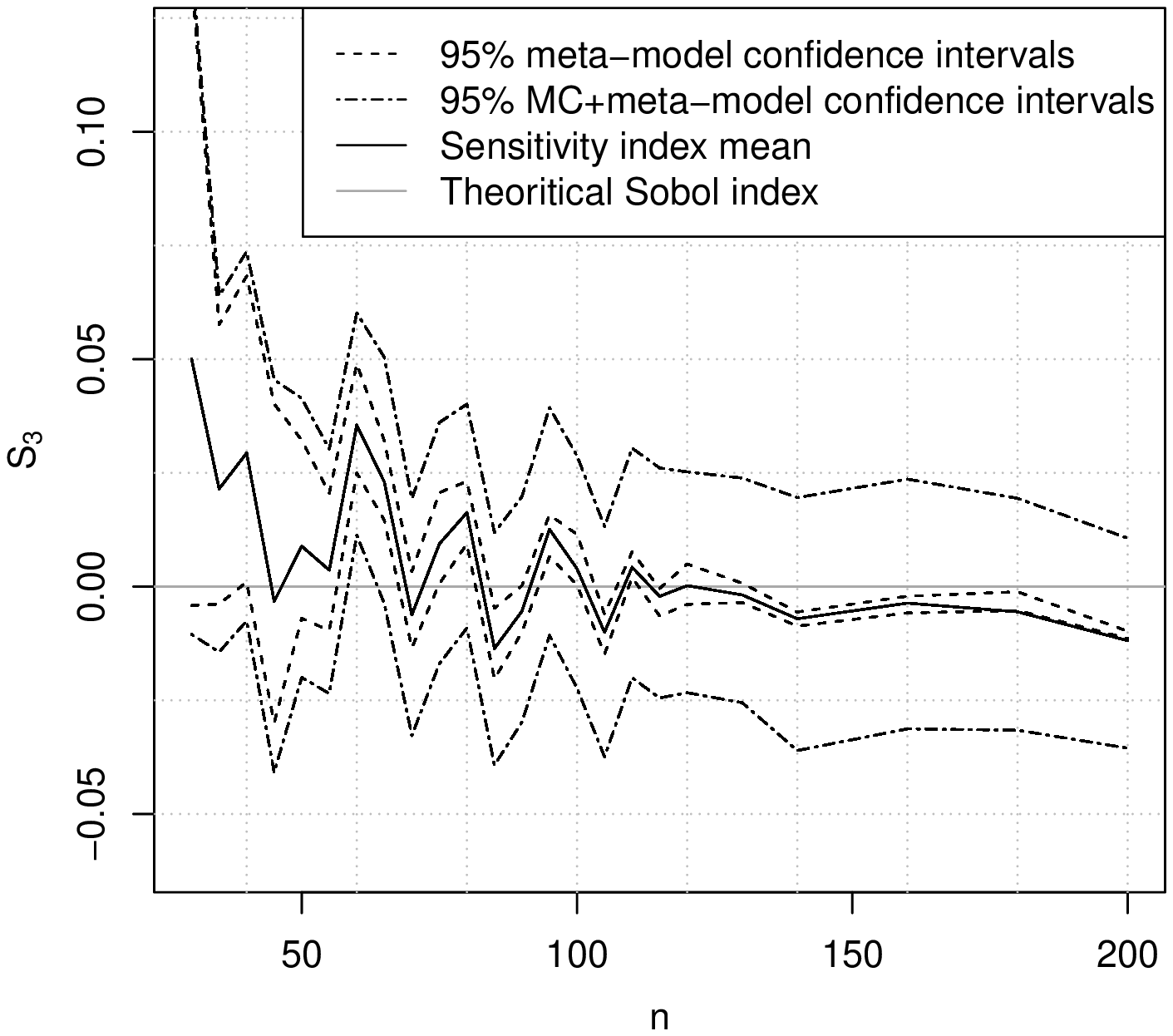}
	\captionsetup{format=hang,justification=centering}
                \caption[hang,center]{}
                \label{c:figS}
        \end{subfigure}
      \caption{Sensitivity index estimates when $n$ increases. The solid lines represent the  means of the sensitivity index estimators. The dashed lines represent the 2.50\% and 97.50\% confidence intervals taking into account only the uncertainty related  to the surrogate modeling. The dashed-dotted lines represent the 2.50\% and 97.50\% confidence intervals taking into account both the uncertainty related  to the surrogate modeling and the one related  to the Monte-Carlo integrations. The horizontal gray lines represent the true values of $S_1$ (\subref{a:figS}), $S_2$ (\subref{b:figS}) and $S_3$ (\subref{c:figS}).}\label{figS}
\end{figure}

Then, we estimate the $0.05$ and $0.95$ quantiles of $\left( \hat{S}^i_{m,n,k,1}\right)$, $k=1,\dots,N_Z$ for each $i=1,2,3$ with  a bootstrap procedure. The resulting quantiles represent the uncertainty related to the surrogate modeling. 
Furthermore, we estimate the $2.50\%$ and $97.50\%$ quantiles of $\left( \hat{S}^i_{m,n,k,l} \right)$, $k=1,\dots,N_Z$, $l=1,\dots,B$  with a bootstrap procedure too.  These quantiles represent the total uncertainty of the index estimator.
Figure \ref{figS} illustrates the result of this procedure for different numbers  of observations $n$. We see in Figure \ref{figS} that for small values of $n$, the error related  to the surrogate modeling dominates. Then, when $n$ increases, this error decreases and it is the one related  to the Monte-Carlo integrations which is the largest. This emphasizes that it is worth to adapt the number of Monte-Carlo particles $m$ to the number of observations $n$. Finally, we highlight that the equilibrium between the two types of uncertainty does not occur for the same $n$  for the three indices. Indeed, it is around $n=100$ for  $S_1$, $n=150$ for  $S_2$ and around $n=75$ for $S_3$. We observe that the smaller  the index is, the larger  its Monte-Carlo estimation error is.

\subsection{Optimal Monte-Carlo resource when $n$ increases}\label{optimalmishi}

We saw in the previous subsection that the equilibrium between the error related to the Monte-Carlo integrations and the one related to the surrogate modeling depends on the considered sensitivity index. The purpose of this subsection is to determine this equilibrium for each index. To perform such analysis, we use the method presented in Subsection \ref{optimalm}.

Let us consider a sample $\left( \hat{S}^i_{m,n,k,l}\right)$, $m=30,\dots,200$,  $k=1,\dots,N_Z$, $l=1,\dots,B$, $i=1,2,3$, generated with Algorithm  \ref{algoSobol1} and using the Monte-Carlo estimator presented in (\ref{Janonestim}). For each pair $(m,n)$ we can evaluate the variance $\hat{\sigma}^2_{Z_n}\left({S}^i_{m,n }\right)$, $i=1,2,3$, related to the meta-modeling with Equation (\ref{uncertaintymodel}) and the variance $\hat{\sigma}^2_{MC}\left({S}^i_{m,n }\right)$, $i=1,2,3$, related to the Monte-Carlo integrations  with Equation (\ref{uncertaintyMC}).
We state that the equilibrium between the two types of uncertainty corresponds to the case
\begin{equation}\label{equilibrium}
\hat{\sigma}^2_{Z_n}\left({S}^i_{m,n }\right) = \hat{\sigma}^2_{MC}\left({S}^i_{m,n }\right).
\end{equation}

\begin{figure}[!ht]
\begin{center}
\includegraphics[width = 8 cm]{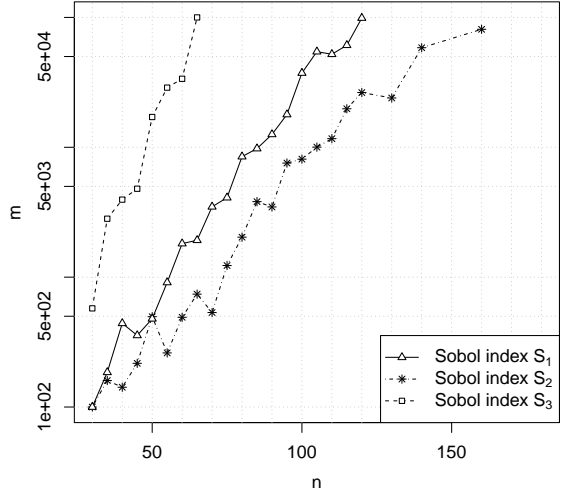}
\caption{Relation between the number of observations $n$ and the number of Monte-Carlo particles $m$ such that the error related to the meta-modeling and the one related  to the Monte-Carlo scheme have the same order of magnitude.}
\label{figEquilibre}
\end{center}
\end{figure}

We present in Figure \ref{figEquilibre} the pairs $(m,n)$ such that the equality (\ref{equilibrium}) is satisfied. We see that the smaller is the sensitivity index, the more important is the number of particles $m$ required to have the equilibrium. Furthermore, we note that the curve increases extremely quickly  for the index $S_3  = 0$. Therefore, it could be  unrealistic to consider the equilibrium for this case, especially when $n$ is important (i.e. $n > 100$).

The presented analysis is of practical interest since it provides  the appropriate number of Monte-Carlo particles $m$ for the sensitivity index estimation in function of the number of observations $n$.  Furthermore, in the framework of computer experiments, the observations are often time-consuming and  $n$ cannot be large.  Therefore, we look for a number of particles $m$ such that the variance $\hat{\sigma}^2_{Z_n}\left({S}^i_{m,n }\right) $ related  to the meta-modeling is smaller than the one of the Monte-Carlo integration $\hat{\sigma}^2_{MC}\left({S}^i_{m,n }\right)$. However, we saw that it could be unfeasible for  some values of  sensitivity index. In this case a compromise must necessarily   be done.

\subsection{Coverage rate of the suggested Sobol index estimator}\label{coverageMC}

\begin{figure}[!ht]
\begin{center}
\includegraphics[width = 8 cm]{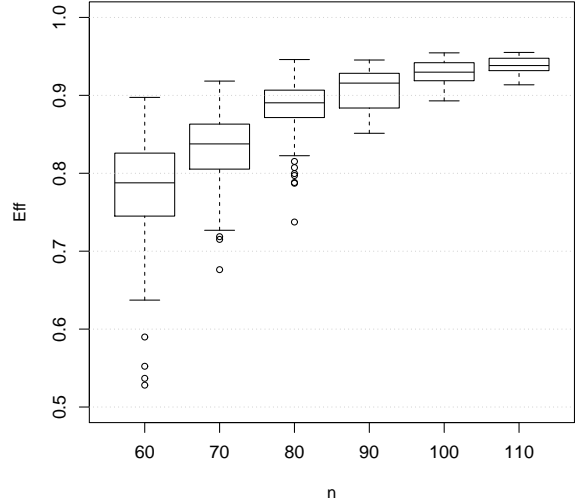}
\caption{Convergence of the model efficiency when the number $n$ of observations increases. For each number of observations $n$, 200 LHS are randomly sampled.}
\label{efficiencylast}
\end{center}
\end{figure}

Algorithm \ref{algoSobol1}  in Subsection \ref{SobolEstimation2} allows for obtaining a sample  $\left( \hat{S}^i_{m,n,k,l}\right)$,    $k=1,\dots,N_Z$, $l=1,\dots,B$ of the estimator of $S_i$ for each $i=1,2,3$. The purpose of this subsection is to verify  the relevance of the confidence intervals provided by $\left( \hat{S}^i_{m,n,k,l}\right)$. To perform such analysis, we  generate $200$  random LHS $(\D_{n,j})_{j=1,\dots,200}$ for different numbers of observations $n$. For each $\D_{n,j}$, we build a kriging model with the procedure presented in Subsection \ref{krigbuilding} and we generate a sample  $\left( \hat{S}^i_{m,n,k,l}\right)$,    $k=1,\dots,N_Z$, $l=1,\dots,B$,  with $B=200$ and $N_Z=300$.
The efficiency of the different kriging models with respect to the number of observation $n$ is presented in Figure \ref{efficiencylast}.
 From this sample, we evaluate the $2.50\%$ and $97.50\%$ quantiles with a bootstrap procedure and we check if the true value of $S_i$ is  covered by  these two quantiles. At the end of the procedure, the ratio between the number of confidence intervals covering the true value  of $S_i$ and the total number of confidence intervals (i.e. 200) has to be close to $95\%$ for each $n$.

Furthermore, to perform the analysis we use different values of $m$ according to the  procedure   presented in Subsection \ref{optimalm} for $S_1$ and $S_2$ (i.e.  such that the variance related  to the meta-modeling  has the same order of magnitude than the one related  to the Monte-Carlo integrations). For $S_3$,  the number of Monte-Carlo particles  $m$ increases too quickly with respect to $n$   to use the method presented in Subsection \ref{optimalm}. Therefore we fix $m$ to the values presented  in Table \ref{optm}. We note that the values of $m$ for $S_3$ are larger than the ones   for $S_1$ and $S_2$.
\begin{table}[H]
\begin{center}
\begin{tabular}{|c||c|c|c|c|c|c|}
\hline
$n$ & 60 & 70 & 80 & 90 &  100 & 110 \\
\hline
$m$ & 1,000 & 3,000 & 5,000 & 10,000 & 40,000 & 60,000  \\
\hline
\end{tabular}
\end{center}
\caption{Numbers of Monte-Carlo particles $m$ for different values of the number of observations $n$ for the estimation  of $S_3$.}
\label{optm}
\end{table}

\begin{figure}[!ht]
        \centering
        \begin{subfigure}[b]{0.45\textwidth}
                \centering
                \includegraphics[width = 6 cm, height = 5cm]{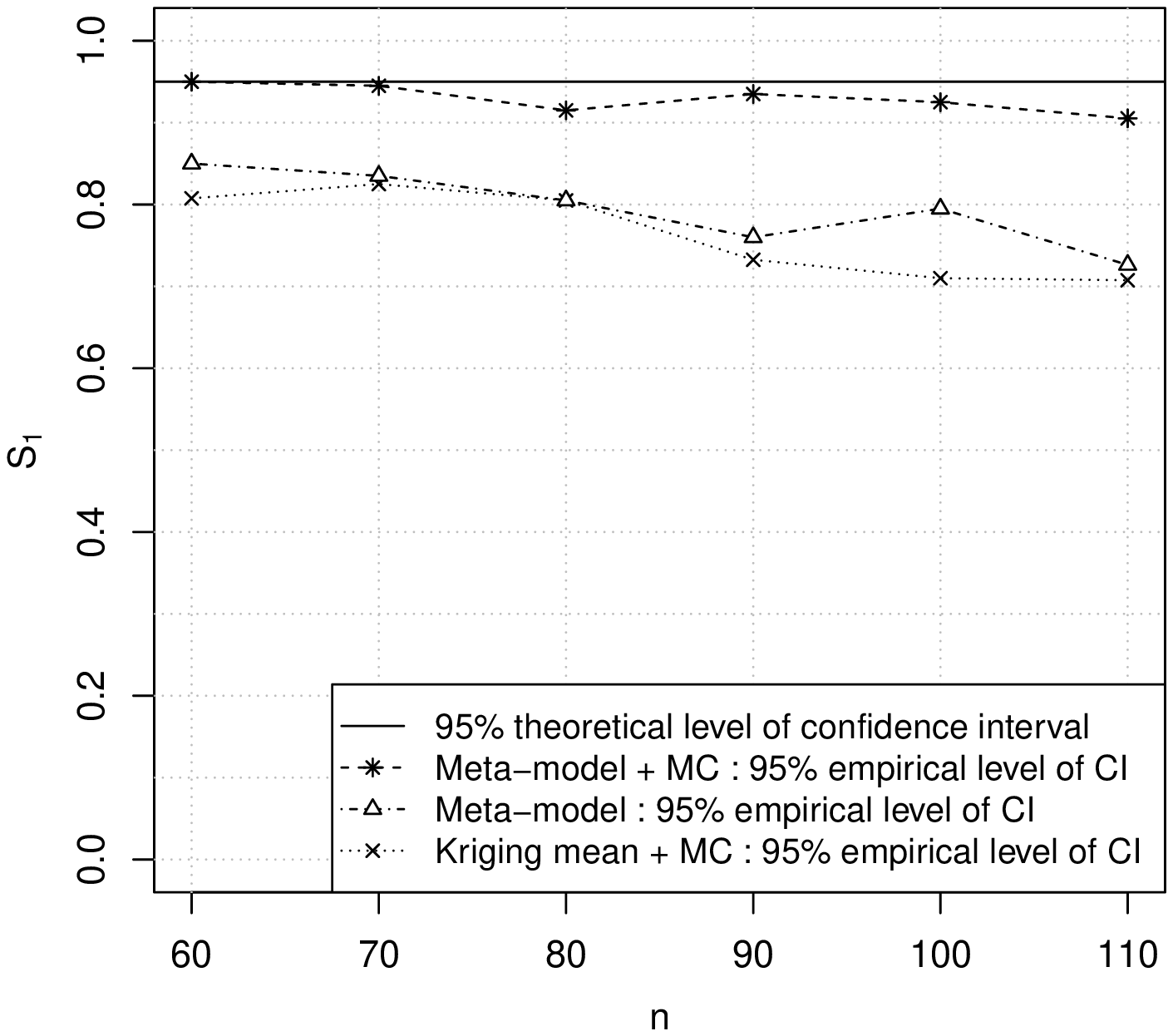}
	\captionsetup{format=hang,justification=centering}
                \caption[hang,center]{}
                \label{a:coverageS}
        \end{subfigure}%
	~
        \begin{subfigure}[b]{0.45\textwidth}
                \centering
                \includegraphics[width = 6 cm, height = 5cm]{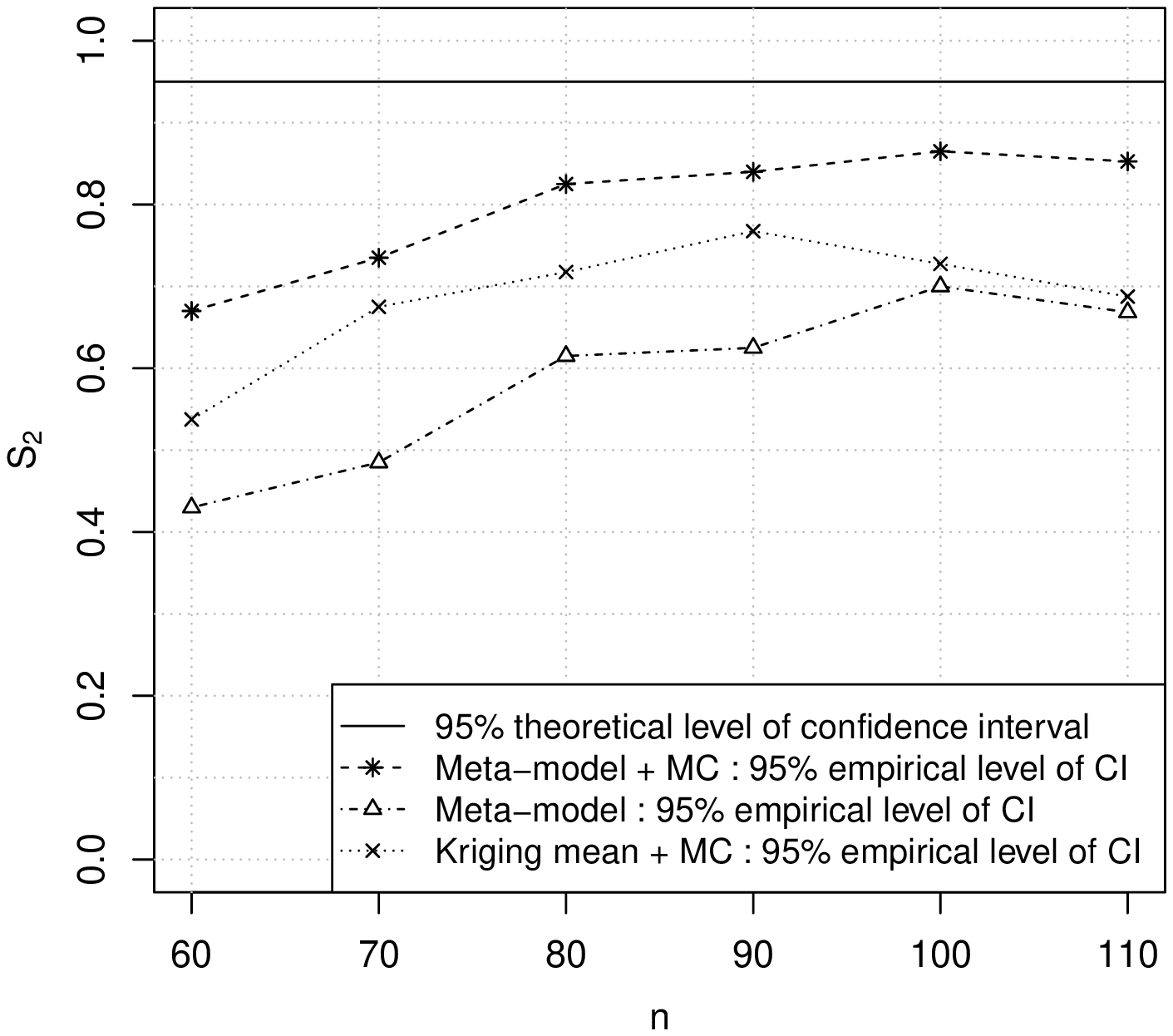}
	\captionsetup{format=hang,justification=centering}
                \caption[hang,center]{}
                \label{b:coverageS}
        \end{subfigure}
 
        \begin{subfigure}[b]{0.45\textwidth}
                \centering               
                \includegraphics[width = 6 cm, height = 5cm]{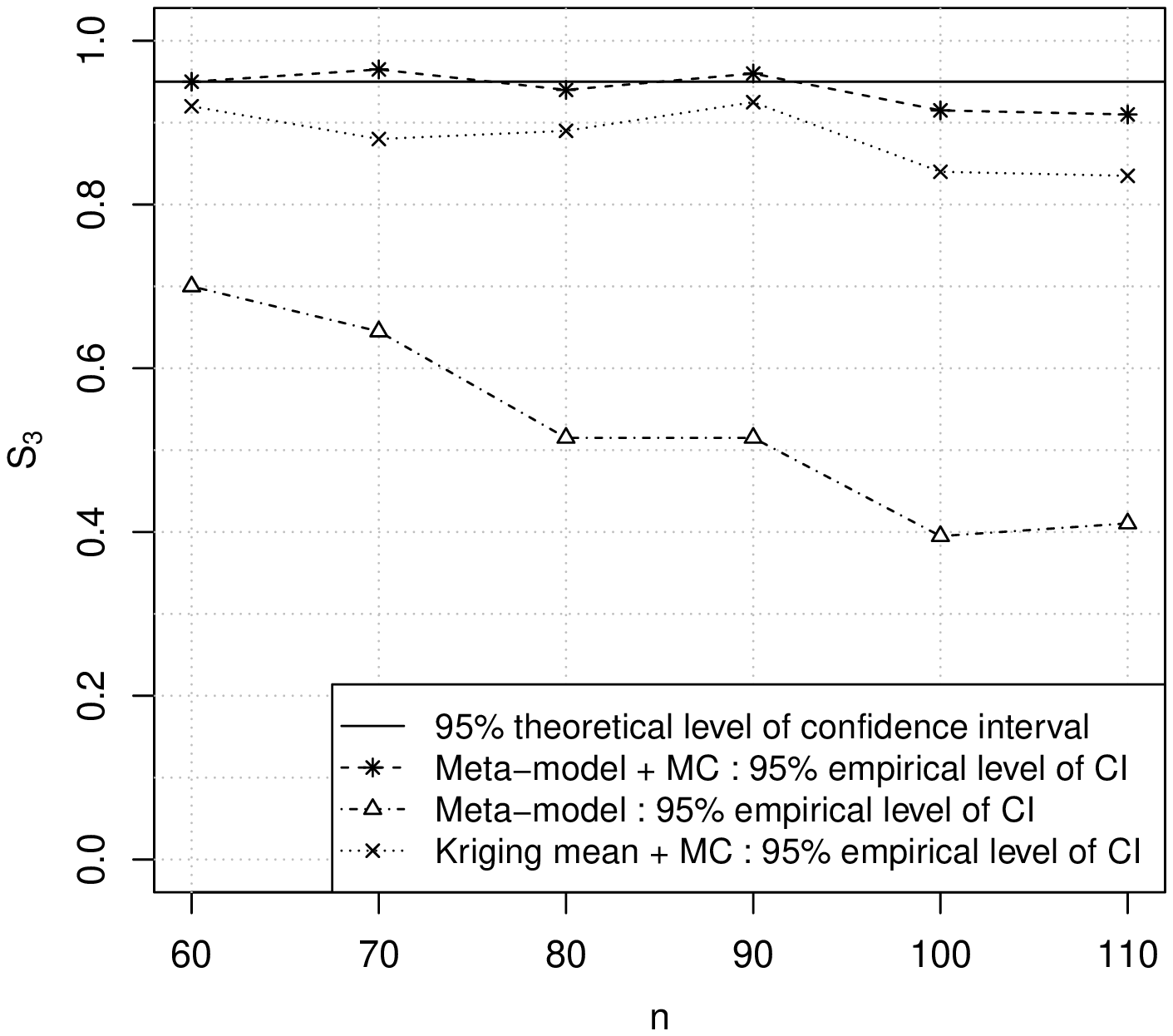}
	\captionsetup{format=hang,justification=centering}
                \caption[hang,center]{}
                \label{c:coverageS}
        \end{subfigure}
      \caption{Empirical 95\% confidence intervals with respect to the number of observations $n$ for $S_1$ (\subref{a:coverageS}), $S_2$ (\subref{b:coverageS}) and $S_3$ (\subref{c:coverageS}).  The empirical coverage rates are evaluated from 200 kriging models built from different random LHS.}\label{coverageS}
\end{figure}

The empirical 95\%-confidence intervals as a function of the number of observations $n$  are presented in Figure \ref{coverageS}. 
We study three cases:
\begin{enumerate}
\item The confidence intervals are built from $\left( \hat{S}^i_{m,n,k,l}\right)$,    $k=1,\dots,N_Z$, $l=1,\dots,B$. Therefore, it takes into account both the uncertainty related to the meta-model and the one related to the Monte-Carlo estimations.
\item  The confidence intervals are built from $\left( \hat{S}^i_{m,n,k,1}\right)$,    $k=1,\dots,N_Z$. In this case,  we  do not use the bootstrap procedure to evaluate the uncertainty due to the Monte-Carlo procedure. Therefore, we only take into account the one due to the meta-model.
\item  The confidence intervals are built from the estimator $\tilde{\mathcal{S}}_{m,n}^{X^{d_1}}$ (\ref{janonGPmean}) with  a bootstrap procedure. Here, we estimate the Sobol indices with  the kriging mean and we do not infer from the uncertainty of the meta-model. Therefore, we only take into account the uncertainty related to the Monte-Carlo estimations.
\end{enumerate}

 We see in Figure \ref{coverageS} that the confidence intervals provided by the approach presented in Section \ref{chap6sec3} are well evaluated for indices $S_1$ and $S_3$.
Furthermore,  they are underestimated when we take into account only the meta-model or the Monte-Carlo   uncertainty. 
This highlights the relevance of the suggested approach to perform uncertainty quantification on the Sobol index estimates.
However, the coverage rate is underestimated for index $S_2$.  This is even worst if we only consider the meta-model error.  This may be due to a poor learning in the   direction  $x_2$ for the the surrogate model. This emphasizes that the suggested method is valid only if the kriging variance well represents the modeling error.

\section{Application of multi-fidelity  sensitivity analysis}\label{sec:application}

In this section, we illustrate the multi-fidelity co-kriging based sensitivity analysis presented in Section \ref{chap6sec4} on an example about a spherical tank under internal pressure.

\subsection{Presentation of the problem}

\begin{figure}[H]
\begin{center}
\vskip-04ex
\includegraphics[trim = 7cm 3cm 1cm 1cm , width = 0.30\linewidth]{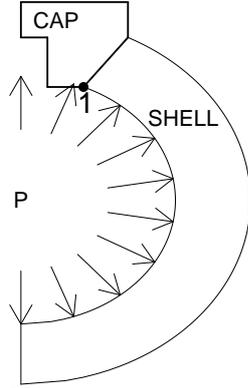}
\vskip-03ex
\caption{Scheme of the  spherical tank under pressure. }
\label{scheme}
\end{center}
\end{figure}

The scheme of the considered tank is presented in Figure \ref{scheme}.  We are interested in the von Mises stress on the  point labeled 1  in Figure \ref{scheme}. It corresponds to the point where the stress is maximal. The  von Mises stresses are of interest since the material yielding occurs when they reach  the critical yield strength.

The system illustrated in Figure \ref{scheme} depends on 8 parameters:
\begin{itemize}
\item $P \, (MPa) \in [30, 50]$: the value of the  internal pressure.
\item $R_{int} \, (mm) \in [1500, 2500]$:  the length of the internal  radius of the shell.
\item $T_{shell} \, (mm) \in [300, 500]$: the thickness of the shell.
\item $T_{cap} \, (mm) \in [100, 300]$: the thickness of the cap.
\item $E_{shell} \, (GPa) \in [63, 77]$: the Young's modulus of the shell material.
\item $E_{cap} \, (GPa) \in [189, 231]$: the Young's modulus of the cap material.
\item $\sigma_{y,shell} \, (MPa) \in [200, 300]$: the yield stress of the shell material.
\item $\sigma_{y,cap} \, (MPa) \in [400, 800]$: the yield stress of the cap material.
\end{itemize}
The von Mises stress $z_2(x)$, $x=(P, R_{int}, T_{shell}, T_{cap}, E_{shell}, E_{cap}, \sigma_{y,shell}, \sigma_{y,cap})$ is provided by an finite elements code, called Aster, modelling the tank under pressure. The material properties of the shell correspond to high quality aluminums and the ones of the cap corresponds to steels from classical to high quality. 

The cheaper version $z_1(x)$ of $z_2(x)$ is obtained by the 1D simplification of the tank corresponding to a perfect spherical tank, i.e. without the cap:
\begin{displaymath}
z_1(x) = \frac{3}{2} \frac{(R_{int}+T_{shell})^3}{(R_{int}+T_{shell})^3 - R_{int}^3}P
\end{displaymath}

\subsection{Multi-fidelity model building}

We present here the construction of the model presented in Section \ref{chap6sec4}. 
%For the implementation, we use the R CRAN package ``MuFiCokriging.1.2'' which computes the equations presented in Section \ref{chap6sec4} and provides maximum likelihood estimates for the hyper-parameters of the covariance kernels.

First, we build two LHS design sets $\tilde{\D}_1$ and $\D_2$ of size $n_1 \times 8$ and $n_2 \times 8$ optimized with respect to the centered $L_2$-discrepancy criterion,  with $n_1 = 100$ and $n_2 =20$.
We note that the input parameter $x$ is normalized so that the measure $\mu(x)$ of the input parameters  is uniform on $  [0,1]^8$.
In order to respect the nested property for the experimental design sets, we remove from $\tilde{\D}_1$ the $n_2$ points  that are the closest to those of $\D_2$ and we set that $\D_1$ is the concatenation of $\D_2$ and $\tilde{\D}_1$. This procedure ensures that $\D_2 \subset \D_1$ without operating any transformation on $\D_2$.

Second, we run the expensive code $z_2(x)$ on points in $\D_2$ and the coarse code $z_1(x)$ on points in $\D_1$. The CPU time of the expensive code is around 1 minute. Furthermore,  in order to have a fair illustration, we consider that the CPU time of the coarse code $z_1(x)$ is not negligible and we restrict its runs to $n_1 = 100$.

Third, we use tensorised $5/2$-Mat\'ern covariance kernels  for $\sigma_1^2 r_1(x,\tx)$ and   $\sigma_2^2 r_2(x,\tx)$ with characteristic length scales $(\theta_1^i)_{i=1, \dots,8}$ and $(\theta_2^i)_{i=1, \dots,8}$. Furthermore, we set that the regression functions are constants, i.e. $\ff_1(x) = 1$ and $\ff_2(x)=1$.

The estimates of the characteristic length scales are given in Table \ref{charscales}.

\begin{table}[H]
\begin{center}
\begin{tabular}{|c||c|c|c|c|c|c|c|c|}
\hline
$\hat{\ttheta}_1$ & 1.71 & 1.38 &  1.97 &  1.98 &  1.98 &  1.99 &  1.95 & 1.41  \\
\hline
$\hat{\ttheta}_2$ & 1.83 &  1.89 & 0.5 &  1.93 &  1.93 &  0.64 &  1.89 &  0.79 \\
\hline
\end{tabular}
\end{center}
\caption{Maximum likelihood estimates of the characteristic length scales of the tensorised $5/2$-Mat\'ern covariance kernels used in the multi-fidelity co-kriging model.  $\hat{\ttheta}_1$ represents the estimates for the code level 1 and $\hat{\ttheta}_2$   represents the ones for the bias between the code levels 1 and 2.}
\label{charscales}
\end{table}

The estimates of the characteristic length scales given in Table \ref{charscales} show that the model is very smooth.  Then, Table \ref{mufiparam} gives the posterior mean of the parameters $(\rho_1, \bbeta_2)$ and $\bbeta_1$ and the restricted maximum likelihood estimate of $\sigma_1^2$ and $\sigma_2^2$.

\begin{table}[H]
\begin{center}
\begin{tabular}{|c||c|}
\hline
$\hat{\bbeta}_1$ & 148.67 \\
\hline
$(\hat{\rho}_1, \hat{\bbeta}_2)$ & (0.92, 57.61) \\
\hline
\end{tabular}
\begin{tabular}{|c||c|}
\hline
$\hat{\sigma}_1^2$ & 495.63 \\
\hline
$\hat{\sigma}_2^2$ & 551.07  \\
\hline
\end{tabular}
\end{center}
\caption{Posterior means of the trend parameters $\bbeta_1$ and $\bbeta_2$ and the adjustment parameter  $\rho_1$ and maximum likelihood estimates of the variance parameters $\sigma_1^2$ and $\sigma_2^2$.}
\label{mufiparam}
\end{table}

The parameter estimates presented in  Table \ref{mufiparam}  show that there is an important bias between the cheap code and the expensive code since $\hat{\bbeta}_2 \approx 58$ whereas the trend of the cheap code is $\hat{\bbeta}_1 \approx 150$. In particular, it is greater than the standard deviation of the bias which is $\hat{\sigma}_2 \approx 23$. 
Then, the  posterior mean of the adjustment parameter $\hat{\rho}_1 = 0.92$ does not indicate a perfect correlation between the two levels of code. Indeed, the estimated correlation between $z_2(x)$ and $z_1(x)$ is $0.77$. Furthermore their estimated variance equals $1514$ for $z_2(x)$ and $810$ for $z_1(x)$. In fact, the adjustment parameter:
\begin{displaymath}
{\rho}_1  = \frac{\mathrm{cov}(Z_2(x), Z_1(x))}{\mathrm{var}(Z_1(x))}
\end{displaymath}
 represents both the correlation degree and the scale factor between the  codes $z_2(x)$ and $z_1(x)$.  

Finally, we can estimate the accuracy of the suggested model with a Leave-One-Out cross validation procedure.  From  the Leave-One-Out errors, we estimate the Nash-Sutcliffe model efficiency $\mathit{Eff}_{LOO} = 83\%$. This means that  the suggested multi-fidelity co-kriging model explains  $83\%$  of the variability of the model. We note that the closer $\mathit{Eff}_{LOO}$ is to 1, the more accurate is the model.  Therefore, we have an excellent model despite the small number of observations $n_2 = 20$ used for the expensive code $z_2(x)$. In order to strengthen this result, we test the multi-fidelity model on an external test set of $7,000$ points and the estimated   efficiency is $86\%$ which is even  better. 

\subsection{Multi-fidelity sensitivity analysis}

Now let us perform a multi-fidelity sensitivity analysis using the approach presented in Subsection \ref{chap6sec40}.  We are interested in the first-order sensitivity indices.

The principle of the method is to sample from the distribution   (\ref{estimatormufi}) using Algorithm \ref{algoSobol2}.  We note that we use  the Monte-Carlo  estimator (\ref{Janonestim}) instead of (\ref{Sobolestim}) since it is asymptotically efficient for the first-order indices
We repeat the algorithm \ref{algoSobol2} to have  $N_Z = 200$ realizations of the predictive distribution $[Z_2(x) | \ZZ^{(2)} = \zz^{(2)}, \ssigma^2]$ and for each realization we generate $B=150$ bootstrap samples. Furthermore,  we choose $m=20,000$ for the Monte-Carlo sampling size  so that the error due to the Monte-Carlo integrations is negligible compared to the one due to the surrogate modelling (see Subsection \ref{optimalm} and \ref{optimalmishi}).

\paragraph{Sensitivity analysis for the cheap code.\\}

First, let us present the result of the sensitivity analysis for the cheap code. As emphasized in Subsection \ref{chap6sec40}, once samples  with respect to the distribution  $[Z_2(x) | \ZZ^{(2)} = \zz^{(2)}, \ssigma^2]$ are available,   samples for  $[Z_1(x) | \ZZ^{(1)} = \zz^{(1)}, \sigma_1^2]$  are also available. Therefore, from them we can perform a sensitivity analysis as presented in  Section \ref{chap6sec3}. Moreover, from  the explicit formula of $z_1(x)$ we expect that only the three variables $P$, $R_{int}$ and $T_{shell}$ have an impact on the output.

\begin{figure}[H]
\begin{center}
\vskip-04ex
\includegraphics[ width = 0.7\linewidth]{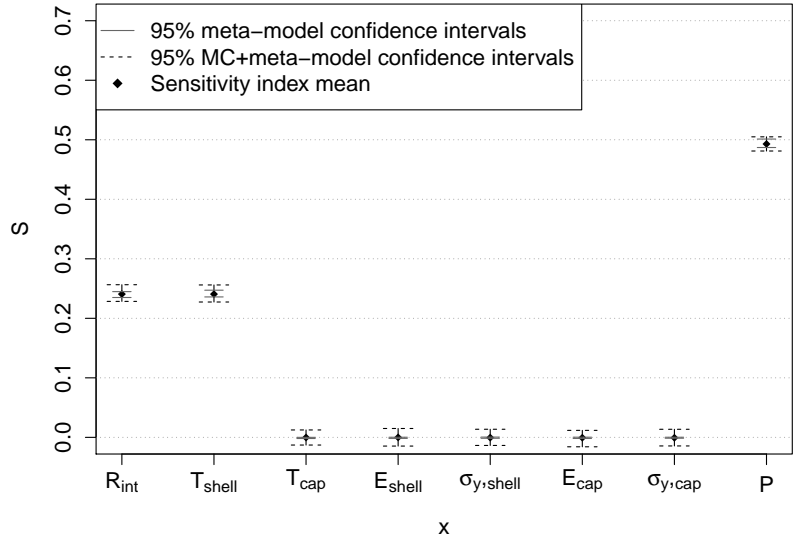}
\vskip-03ex
\caption{Kriging based sensitivity analysis for the cheap code. The diamonds represent the means of the first-order sensitivity index estimators, the solid gray  lines represent  the 95\% confidence intervals taking into account only the meta-modelling uncertainty and the dashed lines represent  the 95\% confidence intervals taking into account  the uncertainty due to both the Monte-Carlo integrations and the meta-modelling.  The means and the confidence intervals are obtained with Algorithm \ref{algoSobol1}. }
\label{ASz1}
\end{center}
\end{figure}

The result of the sensitivity analysis for the cheap code $z_1(x)$ is given in Figure \ref{ASz1}. We see in Figure \ref{ASz1} that only the three parameters $P$, $R_{int}$ and $T_{shell}$ are influent as expected. Furthermore, the internal pressure is the most important parameter whereas the geometrical parameter $R_{int}$ and $T_{shell}$ have equivalent impact on the output. The sum of  the first-order sensitivity index means informs us that $97\%$ of the variability of the output is explained by the first-order indices. The interactions between the parameters are thus negligible. Further, we see that the confidence intervals are tight and that the uncertainty on the Sobol index estimator is essentially due to the Monte-Carlo integrations. This means that the model's error on the cheap code is very low. 

\paragraph{Sensitivity analysis for the expensive code.\\}

\begin{figure}[H]
\begin{center}
\vskip-04ex
\includegraphics[ width = 0.7\linewidth]{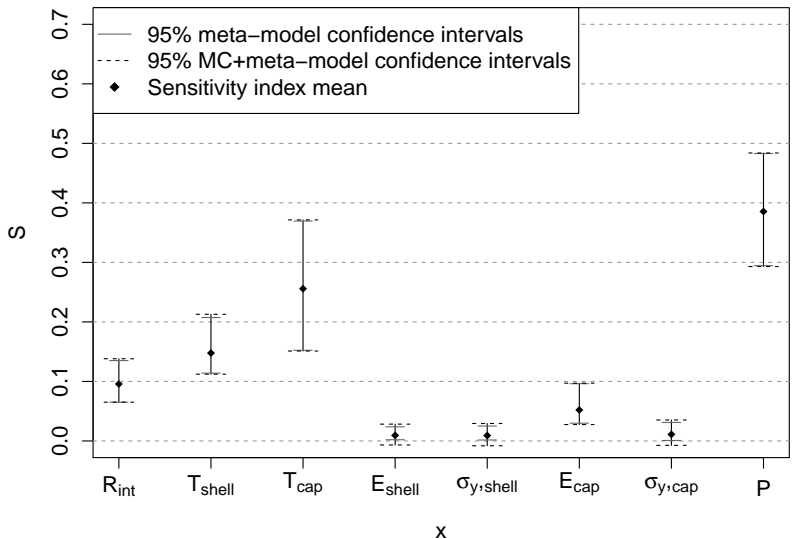}
\vskip-03ex
\caption{Co-kriging based sensitivity analysis for the expensive code. The diamonds represent the means of the first-order sensitivity index estimators, the  solid gray  lines represent  the 95\% confidence intervals taking into account only the meta-modelling uncertainty and the dashed lines represent  the 95\% confidence intervals taking into account  the uncertainty due to both the Monte-Carlo integrations and the meta-modelling.   }
\label{ASz2}
\end{center}
\end{figure}

Second, we perform a sensitivity analysis for the expensive code $z_2(x)$ using the predictive distribution $[Z_2(x) | \ZZ^{(2)} = \zz^{(2)}, \ssigma^2]$. The result of the analysis is presented in Figure \ref{ASz2}.

We see in Figure \ref{ASz2} that the result of the sensitivity analysis for the expensive code is substantially different than the one for the cheap code. First, the importance measure of  the parameters $P$, $R_{int}$ and $T_{shell}$ decreases although the internal pressure $P$ remains the most influent parameter. Second, the material parameters $E_{shell}$, $E_{cap}$, $\sigma_{y,shell}$ and $\sigma_{y,cap}$ have still a negligible influence except for the rigidity of the cap $E_{cap}$. Then, the most noticeable difference is for the thickness of the cap $T_{cap}$ which is now the second most important parameter. Then the sum of the index estimator means equals $96.7\%$. This means that the first-order indices still explain the main part of the model variability.

The hierarchy between the parameters can be easily interpreted. Indeed, the coarse code corresponds to the approximation of the  tank without the cap. Therefore, it is natural that the parameters related to the cap  have no influence. On the contrary, for the expensive code, we are interested in the von Mises stress at the junction between the cap and the shell. Consequently, the parameters related to the cap have now an influence. However, it was difficult   to have a prior on the impact of the cap. We deduce from this analysis that it is in fact  very  important.

Influences of material parameters are negligible because the model stands in the regime of elastic deformations. It is thus physically coherent. In fact, they would be more influent in a plastic deformation regime which can occur for more important internal pressure $P$.

The other important differences between the two sensitivity analysis is the magnitude of the confidence intervals. Indeed, we  see in Figure \ref{ASz2} that, contrary to the cheap code, the confidence intervals for the sensitivity index estimators of    the expensive code are very large. Therefore,  despite the good  multi-fidelity approximation of the expensive code, we have an important uncertainty on it. This is natural since we only use 20 runs of $z_2(x)$  to learn it. 
Finally, we note that the most important uncertainty is for $T_{cap}$. This is explained by the fact that this parameter is not considered by the cheap code. Therefore, $z_1(x)$ brings no information about $T_{cap}$ contrary to $R_{int}$, $T_{shell}$ and $P$.

\section{Conclusion}

This paper deals with the sensitivity analysis of   complex computer codes using Gaussian process regression. The purpose of the paper is to build Sobol index estimators taking  into account both the uncertainty due to the surrogate modelling and the one due to the numerical evaluations of the variances and covariances involved in the Sobol  index definition. The aim is to provide relevant confidence intervals for the index  estimator.

To  provide such estimators, we suggest a method which mixes  a Gaussian process regression model with   Monte-Carlo based integrations. From it,  we can quantify  the impact of both  the Gaussian process regression  and  the Monte-Carlo procedure  on the index estimator variability. In particular, we present a procedure to balance these two sources of uncertainty.
Furthermore, we suggest numerical methods to avoid ill-conditioned problems and to easily handle   the suggested index estimator.

Then, we propose an extension of the suggested approach for multi-fidelity computer codes. These codes have the characteristic to have coarser but computationally cheaper versions. 
They are of practical interest since they allow  for dealing with the problem of very expensive simulations.
To deal with these codes, we use a multivariate Gaussian process regression model called \emph{Multi-fidelity co-kriging}. 

Finally, we perform several numerical tests which confirm the relevance of this new approach.
We illustrate the suggested strategy on an academic example for the univariate case  and with a real application  on a tank under internal pressure for the multi-fidelity analysis.

From this work, two points can naturally be investigated. First, we could improve the uncertainty quantification for the meta-model. Indeed,  in this paper, we do not  take into account the uncertainty due to the estimation of the hyper-parameters of the covariance kernels. This can imply  an underestimation of the predictive variance and thus it can be worth inferring from these parameters. The natural way to perform such analysis is to use a full-Bayesian approach. 
Second, the meta-model considered is built from a fixed experimental design set. Several methods exist to sequentially add new points on the design  in order to perform optimization, to quantify a probability of failure or to improve the accuracy of the meta-model. However, no methods focus on the error reduction of  the sensitivity index estimates. It would be of practical interest to develop sequential design strategies  for a sensitivity analysis purpose.

\section{Acknowledgments}

Part of this work has been backed by French National Research Agency (ANR) through COSINUS program (project COSTA BRAVA noANR-09-COSI-015) and by the CNRS NEEDS program through ASINCRONE project. We
thank Josselin Garnier for several discussions. All the numerical tests have been performed within the R environment, by using the sensitivity, DiceKriging and MuFiCokriging packages.

\bibliographystyle{siam}
\bibliography{biblio}

\end{document}